\documentclass[leqno,11pt]{amsart}
\usepackage{euscript,amssymb,amsbsy}
\usepackage[arrow,matrix,curve]{xy}
\usepackage{a4wide}
\usepackage{graphicx}

\newtheorem{stuff}{Stuff}[section]
\newtheorem{theorem}[stuff]{\sl Theorem}
\newtheorem{proposition}[stuff]{\sl Proposition}
\newtheorem{lemma}[stuff]{\sl Lemma}
\newtheorem{corollary}[stuff]{\sl Corollary}

\newenvironment{definition}{%
\vskip1ex\refstepcounter{stuff}\trivlist \itemindent 0pt
\item[\hskip\labelsep\sl Definition \thestuff.]%
\ignorespaces}{\endtrivlist\vskip1ex}%
\newenvironment{remark}{%
\vskip1ex\refstepcounter{stuff}\trivlist \itemindent 0pt
\item[\hskip\labelsep\sl Remark \thestuff.]%
\ignorespaces}{\endtrivlist\vskip1ex}%

\newtheorem{s-theorem}[sstuff]{\sl Theorem}
\newtheorem{s-proposition}[sstuff]{\sl Proposition}

\let\rar\rightarrow
\let\lar\longrightarrow

\let\hra\hookrightarrow

\let\lmt\longmapsto
\let\xrar\xrightarrow

\font\tenmsa=msam10 %
\newcommand\hdashpiece{%
{\vrule height2.75pt depth-2.35pt width2.3pt \kern1.7pt}}%
\newcommand\hdashpieces{%
{\hdashpiece\hdashpiece\hdashpiece\hdashpiece}}%
\newcommand\dashto{\mathrel{%
\hdashpiece\hdashpiece\kern-0.4pt\hbox{\tenmsa K}}}%
\newcommand\dashar{\mathrel{%
\hdashpieces\kern-0.4pt\hbox{\tenmsa K}}}%

\let\euf\EuScript 
\let\cal\mathcal
\let\mbb\mathbb

\let\mfrak\mathfrak
\let\bsymb\boldsymbol 
\DeclareFontFamily{OT1}{rsfs}{}
\DeclareFontShape{OT1}{rsfs}{n}{it}{<->rsfs10}{}
\DeclareMathAlphabet{\crl}{OT1}{rsfs}{n}{it}

\let\ovl\overline

\let\unbar\underbar

\let\tld\tilde

\let\nit\noindent

\let\disp\displaystyle
\let\srel\stackrel
\let\vphi\varphi

\let\veps\varepsilon

\let\lan\langle
\let\ran\rangle

\newcommand\Hom{\mathop{\sf Hom}\nolimits}
\newcommand\End{\mathop{\sf End}\nolimits}

\newcommand\Pic{\mathop{\rm Pic}\nolimits}
\newcommand\Spec{\mathop{\rm Spec}\nolimits}
\newcommand\uSpec{\mathop{\unbar{\rm Spec}}\nolimits}
\newcommand\Proj{\mathop{\rm Proj}\nolimits}
\newcommand\uProj{\mathop{\unbar{\rm Proj}}\nolimits}
\newcommand\Ker{{\rm Ker}}

\newcommand\invq{{\slash\kern-0.65ex\slash}}

\newcommand\NS{{\rm NS}}

\newcommand\rk{{\rm rank}}

\long\def\InsertFig#1 #2 #3 #4 #5\EndFig{\scalebox{#1}{\hskip #2 mm
$\vbox to #3mm{\vfil\includegraphics{#4}}#5$}}
\long\def\LabelTeX#1 #2 #3\ELTX{\rlap{\kern#1mm\raise#2mm\hbox{#3}}}

\numberwithin{equation}{section}

\let\l\lambda

\let\O\Omega
\let\si\sigma
\let\Si\Sigma

\let\sm\setminus

\newcommand\s{{\rm s}}
\newcommand\sst{{\rm ss}}
\newcommand\us{{\rm us}}

\newcommand\bbA{{\mbb A}}

\newcommand\bbE{{\mbb E}}

\newcommand\bbV{{\mbb V}}
\newcommand\bone{{1\kern-0.57ex\rm l}}

\newcommand\bE{{\bsymb{\rm E}}}
\newcommand\bG{{\bsymb{\rm G}}}
\newcommand\bV{{\bsymb{\rm V}}}
\newcommand\bY{{\bsymb{\rm Y}}}
\newcommand\bGl{{\bsymb{\rm Gl}}}

\newcommand\sfe{{\sf e}}
\newcommand\re{{\rm{\bf e}}}

\newcommand\T{\mathop{\sf T\kern-.2ex}\nolimits}
\newcommand\N{\mathop{\sf N\kern-.1ex}\nolimits}

\newcommand\pr{\mathop{\rm pr}\nolimits}
\newcommand\Sym{\mathop{\rm Sym}\nolimits}
\newcommand\Grs{{\rm Grass}}
\newcommand\codim{{\rm codim}}

\newcommand\Gl{{\rm Gl}}
\newcommand\Ad{{\rm Ad}}
\newcommand\rO{{\rm O}}
\newcommand\rR{{\rm R}}
\newcommand\eO{{{\euf O}}}
\newcommand\eul{{\sf e}}
\newcommand\BG{{\rm BG}}
\newcommand\EG{{\rm EG}}
\newcommand\ch{{\rm ch}}
\renewcommand\det{{\rm det}}

\title{Quotients of affine spaces for actions of reductive groups}
\author{Mihai Halic}
\address{Institut f\"ur Mathematik, Universit\"at Z\"urich,\newline  
\null\hspace{2ex} 
Winterthurerstrasse 190, CH-8057 Z\"urich, Switzerland}
\begin{document}
\maketitle

This article is devoted to the study of invariant quotients of 
an affine space by the action of a reductive group. There are a 
multitude of examples of varieties arising in this way: the most 
at hand are the toric varieties, obtained for representations 
of tori, and the quiver varieties, which include as very particular 
cases the flag varieties of the general linear group. 

The possible linearizations of the action of a linear group on an 
affine space, and more generally on a normal affine variety, are 
given by the characters of the group. A first issue that we address 
is that of the existence of a fan structure on the space of characters, 
whose cones parameterize the possible invariant quotients which appear. 
\smallskip 

\nit{\sl Theorem}\quad 
{\it 
Let $X$ be a normal, affine $G$-variety. The {\rm GIT}-equivalence 
classes in $\cal X^*(G)_{\mbb R}$ corresponding to the $G$-action 
on $X$ are the relative interiors of the cones of a rational, 
polyhedral fan $\Delta^G(X)$, which we call the {\rm GIT}-fan of $X$.
}
\smallskip 

\nit This result can be viewed as a non-abelian generalization 
of the Gelfand-Kapranov-Zelevinskij decomposition for linear 
representations of tori. 
The analogous problem for projective varieties has been 
investigated in \cite{dh,re}. The name of GIT-fan is borrowed 
from \cite{re}, where the author proves a similar result for 
projective varieties. The strategy for proving our result is to 
reduce the problem to the projective case, where we can apply the 
result of Ressayre.  

In the rest of the paper we have tried to answer various questions 
concerning the geometrical properties of the quotients of an affine 
space $\bbV$ on which $G$ acts linearly. 
The first issue in this direction is that of computing the Chow ring 
of the invariant quotients.\smallskip

\nit{\sl Theorem}\quad 
{\it  
Let $\chi\in\cal X^*(G)$ be a character with the property that 
the corresponding semi-stable locus coincides with the 
properly stable one. Then 
$$
A_*\bigl(\bbV\invq_\chi G\bigr)_{\mbb Q}
\cong
\bigl(A_*^G\bigr)_{\mbb Q}
\biggl/
\wp\bigl\lan
[\bbE]_T
\,;\,
\bbE\text{ is a $(T,\chi)$-unstable component of }\bbV
\bigr\ran_{\mbb Q},
\biggr.
$$
where $\wp$ denotes the projection from the $T$- to the 
$G$-equivariant Chow ring of a point. 
If moreover the ground field is that of the complex numbers and 
$\bbV\invq_\chi G$ is projective, then its cohomology ring is 
isomorphic to the Chow ring.
}
\smallskip

This formula has been 
obtained by Ellingsrud and Str\o mme in \cite{es} under the 
additional assumptions that the group $G$ acting on the affine 
space contains its homotheties, and the stability concept 
corresponds to a character with large weight on the subgroup of 
homotheties. We also answer in positive a question raised in 
$loc.\;cit.$, asking whether the Chow ring of the quotients 
is generated by the Chern classes of `natural' vector 
bundles, induced by representations of the group $G$. 

In section \ref{sct:fam} we use the techniques developed by Seshadri 
in \cite{se}, to show that projective varieties of the form 
$\bbV\invq_\chi G$ actually fit together, in families over separated 
and reduced schemes over $\Spec\mbb Z$, and moreover that the absolute 
case considered in the previous sections corresponds to base change. 
If the base is $\Spec\mbb Z$, we obtain flat families of quotients, 
which will eventually allow us to extend in positive characteristic 
certain cohomological properties which hold in characteristic zero. 
In theorem \ref{thm:family} we compute the Chow ring of families 
of such varieties, parameterized by reduced schemes defined over 
algebraically closed fields, and also, in case we work over the 
complex numbers, the cohomology ring of the total space of the 
fibration. 

We conclude the article by computing the Picard group and the 
ample cone of the geometric quotients that we obtain. Using the 
Hochster-Roberts theorem, we further prove the vanishing of 
higher cohomology groups of nef line bundles. 
\smallskip 

\nit{\sl Theorem}\quad 
{\it  
Let $\chi\in\cal X^*(G)$ be a character such that the codimension 
of the unstable locus ${\rm codim\,}_\bbV\bbV^\us(G,\chi)\geq 2$, 
and $\bbV^\sst(G,\chi)=\bbV^\s_{(0)}(G,\chi)$. There is a prime 
$p(\chi)$ depending only on the {\rm GIT}-chamber of $\chi$, having 
the following property: if the characteristic of the ground field 
is zero or larger than $p(\chi)$, then for any nef, invertible   
sheaf $\crl L\rar\bbV\invq_\chi G$ holds: 
$$
H^i\bigl(\bbV\invq_\chi G,\crl L\bigr)=0,
\;\forall\,i>0,
\quad\text{and}\quad
H^i\bigl(\bbV\invq_\chi G,\crl L^{-1}\bigr)=0,
\;\forall\,i\neq\kappa(\crl L),
$$
where $\kappa(\crl L)$ denotes the Kodaira-Iitaka dimension. 
}
\smallskip

In characteristic zero, we prove in theorem \ref{thm:X} that a 
similar vanishing result holds for the quotients of smooth, affine 
varieties.
\smallskip 

\nit{\small\it Acknowledgements}\quad This work has been partially 
supported by EAGER\ -\ European Algebraic Geometry Research 
Training Network, contract No. HPRN-CT-2000-00099 (BBW). I am very 
indebted to M.~Brodmann for patiently explaining me properties of 
the local cohomology modules, that I needed in section 
\ref{coh-line-bdls}. I thank also S.~Stupariu for useful 
discussions, which helped me to improve the presentation in 
section \ref{sct:fan}.


\section{Stability for linear actions of reductive groups \sf I}
{\label{sct:stab}}

Let $K$ be an algebraically closed field of arbitrary characteristic, 
$G$ a reductive group over $K$, and $G_m:=\Spec K[t,t^{-1}]$ the 
multiplicative group of $K$. We denote $Z(G)^\circ$ the connected 
component of the centre of $G$; we fix a maximal torus $T$ of $G$ 
and denote by $W:=N_G(T)/T$ the corresponding Weyl group. Further, 
we consider a finite dimensional $G$-module $V$, such that the 
representation $\rho:G\rar\Gl(V)$ has finite kernel. 

Let $\bbV:=\Spec\bigl(\Sym^\bullet V\bigr)$ be the affine space 
corresponding to $V$, and consider  
$$
\Si : G\times \bbV\lar \bbV
$$
the natural $G$-action on it (see \cite[section 1.1]{mfk}). 
Since $\Pic(\bbV)=\{\crl O_{\bbV}\}$, the possible linearizations 
of $\Si$ are parameterized by the group of characters $\cal X^*(G)$ 
of $G$: namely, $\chi\in\cal X^*(G)$ defines the linearization  
\begin{align}{\label{si-chi}} 
\Si_\chi:G\times\bigl(\bbV\times\bbA^1_K\bigr)
\lar
\bbV\times\bbA^1_K,
\quad
\Si_\chi\bigl(g,(x,z)\bigr):=\bigl(\Si(g,x),\chi(g)z\bigr). 
\end{align}
In this formula, and overall in the article, we adopt the convention 
that for defining the linearization of a $G$-action in a sheaf we 
indicate the corresponding $G$-action on the geometric line bundle 
defined by the sheaf. In the case above, the geometric line bundle is 
$\eO_\bbV:=\uSpec\bigl(\Sym^\bullet\!\crl O_\bbV^\vee\bigr)
=\bbV\times\bbA^1_K$. A function $f\in K[\bbV]$ is $\Si_\chi$-invariant 
precisely when $f\bigl(\Si(g,x)\bigr)=\chi(g)f(x)$, for all 
$(g,x)\in G\times\bbV$. For $n\geq 0$, we denote by 
$K[\bbV]^G_{\chi^n}$ the vector space of $\Si_{\chi^n}$-invariant 
functions on $\bbV$, and define the algebra 
$$
K[\bbV]^{G,\chi}:=\bigoplus_{n\geq 0}K[\bbV]^G_{\chi^n}.
$$ 
The goal of this section is to describe the (semi-)stable 
locus $\bbV^\sst(G,\chi)$ corresponding to $\Si_\chi$, for 
$\chi\in\cal X^*(G)$. Proposition \ref{prop:sst} below is stated 
explicitely in \cite[proposition 2.5]{ki}, but the proof is very 
sketchy. Since for our computations we need a clear description 
of the unstable locus, we found useful to give a thorough proof 
of the numerical criterion for linear actions. Another remark is 
that in $loc.\;cit.$ the author allows the representation $\rho$ 
to have kernel; however, one notices that, if the 
$\Si_\chi$-semi-stable locus is not empty, then the restriction 
to the connected component of the identity 
$\chi|_{(\Ker\rho)^\circ}\equiv 1$, by trivial reasons. Therefore 
the assumption that the kernel of $\rho$ is finite is not restrictive. 

A one-parameter subgroup (1-PS for short) $\l\in\cal X_*(G)$ 
of $G$ decomposes $V$ into the direct sum of the weight spaces 
$V=\oplus_jV_j$, where $\l(t)|_{V_j}=t^j{\rm Id}_{V_j}$. We observe 
that, except for a finite number of $j$'s, all the $V_j$'s are zero. 
This decomposition breaks $\bbV$ into the direct product 
$$
\bbV=\Spec\bigl(\Sym^\bullet V\bigr)
=\Spec\bigl(\otimes_j\Sym^\bullet V_j\bigr)
=\times_j\mbb V_j,
$$
and we remark that $\l$ acts on $\bbV_j$ by multiplication through 
$t^j$. We denote $(x_j)_j$ the components of a point $x\in\bbV$ 
with respect to this decomposition, and we define
$$
m(x,\l):=
\biggl\{
\begin{array}{ll}
\min\{j\mid x_j\neq 0\}&\text{if }x\neq 0,
\\[1ex] 
0&\text{if }x=0.
\end{array}
\biggr.
$$
\begin{lemma}{\label{lm:us}}
{\rm (i)} For $\chi\in\cal X^*(G)$, the corresponding 
$\chi$-unstable locus is 
$$
\bbV^\us(G,\chi)
=\bigcup_{\hbox{\scriptsize
$\begin{array}{l}
\l\in\cal X_*(T),
\\ 
\lan\chi,\l\ran<0
\end{array}$}}
\kern-3ex
G\cdot\bbE(\l)
=G\cdot\bbV^\us(T,\chi),
$$
where $\bbE(\l):=\{x\in\bbV\mid m(x,\l)\geq 0\}$. It is invariant 
under the action of $N_{\Gl(V)}G$, the normalizer of $G$ in $\Gl(V)$.

\nit{\rm (ii)} For each $\l\in\{\lan\chi,\cdot\ran<0\}$, $\bbE(\l)$ 
is a linear subspace of $\bbV$, which is stable under the parabolic 
subgroup 
$$
P(\l):=\bigl\{
g\in G \mid \lim_{t\rar 0}\l(t)g\l(t)^{-1}\text{ exists in }G
\bigr\}
$$ 
of $G$. As a consequence, $G\cdot\bbE(\l)$ is closed in $\bbV$.

\nit{\rm (iii)} There is a finite set 
$\euf F(\chi)\subset\{\lan\chi,\cdot\ran<0\}$, such that 
$$
\bigcup_{\hbox{\scriptsize
$\begin{array}{l}
\l\in\cal X_*(T),
\\ 
\lan\chi,\l\ran<0
\end{array}$}}
\kern-3ex
G\cdot\bbE(\l)
=
\bigcup_{\l\in\euf F(\chi)}G\cdot\bbE(\l).
$$
It corresponds to the irreducible components of the union on 
the left-hand-side above. 
\end{lemma}

\begin{proof}
We define the representation 
$$
\rho':G\lar \Gl(K\oplus V),\quad
\rho'(g):={\rm diag}\bigl(\chi^{-1}(g),\rho(g)\bigr),
$$
and consider the induced action $\Si':G\times(\bbA^1_K\times\bbV)
\rar \bbA^1_K\times\bbV$, $\Si'\bigl(g,(z,x)\bigr)
=\bigl(\chi(g)z,\Si(g,x)\bigr)$. The corresponding ring of 
invariants is obviously
$$
K[\bbA^1_K\times\bbV]^G=
\biggl\{
F=\sum_{n\geq 0}t^nf_n\,\biggl|\, f_n\in K[\bbV]^G_{\chi^n}
\biggr\}.
$$
We deduce that 
$$
\begin{array}{ll}
x\in\bbV^\us(G,\chi) 
& \Leftrightarrow\; 
f(x)=0,\ \forall\,f\in K[\bbV]^G_{\chi^n}, n\geq 1
\\[1ex]
&\Leftrightarrow\; 
F(1,x)=F(0,x)=f_0(x),\ \forall\,F\in K[\bbA^1_K\times\bbV]^G
\\[1ex]
&\Leftrightarrow\; 
F(1,x)=0,\ \forall\,F\in ({\euf I_{\{0\}\times\bbV}})^G
\\[1ex]
&\Leftrightarrow\; 
\overline{\Si'\bigl(G\times (1,x)\bigr)}\cap 
\bigl(\{0\}\times\bbV\bigr)\neq \emptyset.
\end{array}
$$
We deduce from \cite[theorem 1.4]{ke} that the last condition is 
equivalent to the existence of $\l\in\cal X_*(G)$ such that 
$\lim_{t\rar 0}\l(t)\cdot(1,x)\in \{0\}\times\bbV$. This happens 
precisely when $\langle\chi,\l\ran <0$ and $m(x,\l)\geq 0$. The 
conclusion follows now, because any such $\l$ is contained in a 
maximal torus of $G$, and moreover all maximal tori of $G$ are 
conjugated to $T$. 

(ii) The first part is obvious and we prove the second one: we notice 
that $\bbE(\l)=\{x\in\bbV\mid\lim_{t\rar 0}\l(t)\cdot x
\text{ exists in }\bbV\}$. For $g\in P(\l)$, we may write
$$
\Si\bigl(\l(t),\Si(g, x)\bigr)
=\Si\bigl(\l(t)g\l(t)^{-1},\Si\bigl(\l(t), x\bigr)\bigr),
$$
and deduce that the limit at zero exists, hence $\Si(g,x)\in\bbE(\l)$. 
The last statement is proved in proposition \ref{prop:image}. 

(iii) The number of one parameter subgroups of a torus is countable, 
and a countable union of closed subvarieties of $\bbV$, which is itself 
closed, is actually a finite union. 
\end{proof}

We are now  in position to state the Hilbert-Mumford criterion 
corresponding to our linear problem.

\begin{proposition}{\label{prop:sst}} Fix a character 
$\chi\in\cal X^*(G)$. Then:

\nit{\rm (i)}\quad 
$x\in\bbV^\sst(G,\chi)
\;\Longleftrightarrow\;
\bigl[
\forall\,\l\in\{\lan\chi,\cdot\ran<0\}\,\Rightarrow\,m(x,\l)<0
\bigr].$

\nit{\rm (ii)}\quad 
$x\in\bbV^\s_{(0)}(G,\chi)
\;\Longleftrightarrow\;
\bigl[
\forall\,\l\in\{\lan\chi,\cdot\ran\leq 0\}\,\Rightarrow\,m(x,\l)<0
\bigr].$
\end{proposition}

\begin{proof}
The first equivalence is a straightforward consequence of the 
previous proposition. We are going to prove the second one. 

\nit($\Rightarrow$)\quad Suppose that there exists 
$x\in\bbV^\s_{(0)}(G,\chi)$, which does not fulfill the condition 
on the right-hand-side. Since $x$ is stable, it is automatically 
semi-stable; this fact, together with our hypothesis, implies the 
existence of a 1-PS $\l$ such that $\lan\chi,\l\ran=0$ and 
$m(x,\l)\geq 0$. Then $x_0:=\lim_{t\rar 0}\Si(\l(t), x)$ exists in 
$\bbV$, and differs from $x$, since the stabilizer of $x$ is finite. 

We claim that $x_0\in\bbV^\sst(G,\chi)$: indeed, there is a function 
$f\in K[\bbV]^G_{\chi^n}$, $n\geq 1$, such that $f(x)\neq 0$. Then 
$$
f(x_0)
=\lim_{t\rar 0}f\bigl(\Si(\l(t), x)\bigr)
=\chi^n(\l(t))\cdot f(x)=t^{n\lan\chi,\l\ran}f(x)
=f(x)\neq 0.
$$
But this means that the $G$-orbit of $x$ is not closed in 
$\bbV^\sst(G,\chi)$, which is a contradiction.

\nit($\Leftarrow$)\quad Consider $x\in\bbV$ which satisfies the relation 
on the right-hand-side: then it is automatically semi-stable. We must 
prove that its orbit is closed in $\bbV^\sst(G,\chi)$ and that its 
stabilizer is finite. 

Assume that the $G$-orbit of $x$ is not closed in $\bbV^\sst(G,\chi)$: 
then we find a 1-PS $\l$ such that 
$x_0:=\lim_{t\rar 0}\Si\bigl(\l(t), x\bigr)\in\bbV^\sst(G,\chi)\sm
\bbV^\s_{(0)}(G,\chi)$. 
We deduce that $m(x,\l)\geq 0$, and therefore $\lan\chi,\l\ran>0$. 
For any $f\in K[\bbV]^G_{\chi^n}$, $n\geq 1$, we have
$$
f(x_0)=\lim_{t\rar 0}f\bigl(\Si(\l(t), x)\bigr)
=\lim_{t\rar 0}\chi^n(\l(t))\cdot f(x)
=\lim_{t\rar 0}t^{n\lan\chi,\l\ran}f(x)=0,
$$
which contradicts the fact that $x_0\in\bbV^\sst(G,\chi)$.

Since $Gx:=\Si(G,x)$ is closed in $\bbV^\sst(G,\chi)$, it is closed 
in some $G$-invariant affine subset containing $x$; it follows that 
$Gx$ is affine, and therefore the stabilizer of $x$ is reductive, by 
Matsushima's criterion (see \cite{ri}). Assume now that there is a 
1-PS $\l$ contained in the stabilizer of $x$: then we may chose it 
such that $\lan\chi,\l\ran\geq0$. Obviously $m(x,\l)=0$, which 
contradicts our hypothesis. 
\end{proof}


\section{Stability for linear actions of reductive groups \sf II}
{\label{sct:stab2}}

In this section we will view the affine space $\bbV$ as a 
$G$-invariant, Zariski open subset of the projective space 
$\mbb P(K\oplus V):=\Proj\bigl(\Sym^\bullet(K\oplus V)\bigr)$, 
for an appropriate $G$-action on this latter one. For any character 
$\chi\in\cal X^*(G)$, we prove that the $G$-action on 
$\mbb P(K\oplus V)$ can be linearized in such a way that the 
intersection of the corresponding semi-stable locus with $\bbV$ 
coincides with the $\chi$-semi-stable locus of $\bbV$. The goal of 
this approach is that of allowing us to reduce the issue of chamber 
structure on $\cal X^*(G)$ in the affine case to known results in 
the projective case. 

Let us recall that for any character $\chi\in\cal X^*(G)$, the 
invariant quotient $\bbV\invq_\chi G$ corresponding to the 
linearization $\Si_\chi$ can be described as
$$
\bbV\invq_\chi G=\bbV^\sst(G,\chi)/G
=\Proj\bigl(K[\bbV]^{G,\chi}\bigr),
$$ 
and this one is projective over $\Spec\bigl(K[\bbV]^G\bigr)$. 

Let us denote $\pi:\bbV=\Spec K[\bbV]\rar\Spec K[\bbV]^G=:\bbV/G$ the 
natural projection, and define $\hat 0:=\pi(0)\in\bbV/G$. Further, 
$K[\bbV]^{G,\chi}$ being a subalgebra of $K[\bbV]$ inherits an 
additional graduation given by the degree in $K[\bbV]$ (it becomes a 
bi-graded algebra); we denote by $K[\bbV]^G_{\chi^n,(p)}\subset 
K[\bbV]^G_{\chi^n}$ the submodule homogeneous elements of degree $p$. 

\begin{lemma}{\label{lm:grow}}
For a character $\chi\in\cal X^*(G)$, the following statements hold:

\nit{\rm (i)} There is an integer $c>0$ such that 
$K[\bbV]^G_{\chi^{nc}}=\bigl(K[\bbV]^G_{\chi^c}\bigr)^n$ for all 
$n\geq 1$. 

\nit Assume that $\chi\in\cal X^*(G)$ is such that the condition 
{\rm (i)} holds with $c=1$. 

\nit${\rm (ii_a)}$ In the case $K[\bbV]^G\neq K$, there are integers 
$D,D_\chi\geq 1$ such that the following two conditions are fulfilled:
$$
\biggl\{
\begin{array}{l}
\forall n\geq 1\,\forall x\in\bbV\sm\pi^{-1}(\hat 0)\ 
\exists f\in K[\bbV]^G_{(nD)}\text{ with }f(x)\neq 0;
\\[0.7ex] 
\forall n>D_\chi\,\forall x\in\bbV^\sst(G,\chi)\sm\pi^{-1}(\hat 0)\ 
\exists f\in K[\bbV]^G_{\chi^D,(nD)}\text{ with }f(x)\neq 0. 
\end{array}
\biggr.
$$
Conversely: consider $x\in\bbV$ with the property that there are 
$N>0$ and $f\in K[\bbV]^G_{\chi^N,(Nd)}$, with $d>D_\chi$, such that 
$f(x)\neq 0$. Then $x\in\bbV^\sst(G,\chi)\sm\pi^{-1}(\hat 0)$.

\nit${\rm (ii_b)}$ In the case $K[\bbV]^G=K$, 
$$
\max\bigl\{\deg f\mid f\in K[\bbV]^G_{\chi^n}\bigr\}\leq 
n\cdot\max\bigl\{\deg f\mid f\in K[\bbV]^G_\chi\bigr\}=:nD_\chi, 
\text{ for all }n\geq 1.
$$
\end{lemma}

\begin{proof}
The first statement is a particular case of \cite[lemme 2.1.6]{ega2}, 
applied to the noetherian ring $K[\bbV]^{G,\chi}$. The statement 
$\rm (ii_b)$ is obvious, since in this case $K[\bbV]^G_\chi$, which 
is a finite dimensional $K$-module, generates $K[\bbV]^{G,\chi}$ as 
a $K$-algebra. 

We prove now $\rm (ii_a)$: we consider finite sets of homogeneous 
polynomials $\cal S_0\subset K[\bbV]^G\sm K$ and $\cal S_1\subset 
K[\bbV]^G_\chi$ which generate respectively $K[\bbV]^G$ as a 
$K$-algebra, and $K[\bbV]^{G,\chi}$ as a $K[\bbV]^G$-algebra. 
We define $D$ to be the smallest common multiple of $\deg f$, for 
$f\in\cal S_0$, and $D_\chi:=\max\{\deg f_1\mid f_1\in\cal S_1\}$. 

Then, for any $x\in\bbV\sm\pi^{-1}(\hat 0)$, there is $f\in\cal S_0$ 
such that $f(x)\neq 0$. Raising it to a suitable power, we obtain 
the desired polynomial. 

Consider an integer $n\geq 1$. For $x\in\bbV^\sst(G,\chi)\sm 
\pi^{-1}(\hat 0)$, we find $f_1\in\cal S_1$ such that $f(x)\neq 0$. 
Moreover, the previous step shows that there is a homogeneous 
polynomial $f_0\in K[\bbV]^G$ of degree 
$\deg f_0=(n+D_\chi-\deg f_0)D$ with $f_0(x)\neq 0$. 

We obtain that 
$$ 
f_0f_1^D\in K[\bbV]^G_{\chi^D},\ 
\deg(f_0f_1^D)=\bigl(n+D_\chi\bigr)D,
\text{ and }f_0f_1^D(x)\neq 0. 
$$
Conversely, assume that we have $x\in\bbV$ with the stated property. 
Then, clearly, $x\in\bbV^\sst(G,\chi)$. We claim that the elements 
of $K[\bbV]^G\sm K$ does not vanish at $x$. Assume the contrary, and 
consider a homogeneous decomposition $f=f'+f''$ with 
$f'\in K[s;s\in\cal S_1]$ and $f''\in \cal I[s;s\in\cal S_1]$, 
for $\cal I:=\lan s_0;s_0\in\cal S_0\ran\subset K[\bbV]^G$. Since 
$f(x)\neq 0$, $f'(x)\neq 0$ too. Thus we have found a non-zero 
$f'\in K[\bbV]^G_{\chi^N,(Nd)}$ which can be expressed as a polynomial 
of degree $N$ (because of the $\chi^N$-homogeneity) in the elements 
of $\cal S_1$. But, in this case 
$$
Nd=\deg f'\leq N\cdot\max\{\deg f_1\mid f_1\in\cal S_1\}=ND_\chi,
$$
which contradicts the choice of $d>D_\chi$.
\end{proof}

In what follows we view $\mbb P(K\oplus V)=
\Proj\bigl(\Sym^\bullet(K\oplus V)\bigr)\cong 
\bigl((\bbA^1_K\times\bbV)\sm\{0\}\bigr)\bigl/G_m$. 
We consider the representation $G\rar\Gl(K\oplus V)$, 
$g\lmt{\rm diag}\bigl(1,\rho(g)\bigr)$, and observe that 
for the induced action
\begin{align}{\label{eqn:action}}
\bar\Si:G\times\mbb P(K\oplus V)\lar \mbb P(K\oplus V),
\quad
\bar\Si(g,[a,x]):=[a,\Si(g,x)],
\end{align}
the inclusion $\bbV\subset\mbb P(K\oplus V)$ is $G$-invariant. 
What we wish is to linearize the $G$-action in 
$\crl O_{\mbb P(K\oplus V)}(d)$, for some $d>0$, such that the 
corresponding (semi-)stable points in $\bbV$ are precisely the 
$\Si_\chi$-(semi-)stable points of $\bbV$.

We remark that for an integer $d$, the geometric line bundle 
defined by $\crl O_{\mbb P(K\oplus V)}(d)$ can be described as   
$$
\eO_{\mbb P(K\oplus V)}(d) :=
\uSpec\bigl(\Sym^\bullet\!\crl O_{\mbb P(K\oplus V)}(-d)\bigr)
\cong  
\raise.5ex\hbox{
$\bigl((\bbA^1_K\times\bbV)\sm\{0\}\bigr)\times\bbA^1_K$}
\bigr/
\lower.5ex\hbox{$G_m$},
$$
where $G_m$ acts on $(\bbA^1_K\times\bbV)\times\bbA^1_K$ by 
$t\times(y,z):=(ty,t^dz)\;\;\forall t\in G_m$. 
For $\chi\in\cal X^*(G)$, we linearize the $G$-action on
$\mbb P(K\oplus V)$ in $\crl O_{\mbb P(K\oplus V)}(d)$ by
\begin{align}{\label{eqn:lin}}
\Si_{d,\chi}:G\times \eO_{\mbb P(K\oplus V)}(d)\lar
\eO_{\mbb P(K\oplus V)}(d),
\quad
\bigl(g,\bigl[(a,x),z \bigr]\bigr)\lmt 
\bigl[(a,\Si(g,x)),\chi(g)z \bigr].
\end{align}
In particular, we deduce that there is a natural isomorphism
\begin{align}{\label{split}}
\Pic^G\bigl(\mbb P(K\oplus V)\bigr)
\cong
\Pic\bigl(\mbb P(K\oplus V)\bigr)\oplus\cal X^*(G)
\cong\mbb Z\oplus\cal X^*(G).
\end{align}

\begin{lemma}{\label{lm:lin}}
{\rm (i)} The restriction of the action $\Si_{d,\chi}$ 
to $\eO_{\mbb P(K\oplus V)}(d)|_\bbV$ induces the action 
$\Si_\chi$ on $\eO_\bbV$.

\nit{\rm (ii)} The restriction homomorphism
$\Pic^G\bigl(\mbb P(K\oplus V)\bigr)\lar \Pic^G(\bbV)$
is surjective and corresponds to the projection 
$\mbb Z\oplus\cal X^*(G)\rar\cal X^*(G)$.
\end{lemma}

\begin{proof}
We consider the section $s$ in $\eO_{\mbb P(K\oplus V)}(d)$
given by $s[a,x]:=[(a,x),a^d]$; it defines an isomorphism
$\eO_\bbV\srel{\cong}{\lar}\eO_{\mbb P(K\oplus V)}(d)|_\bbV$,
and we observe that
$\Si_{d,\chi}\bigl(g, s[1,x]\bigr)=\chi(g)\cdot s[1,\Si(g,x)]$ 
for $x\in\bbV$ and $g\in G$.
\end{proof}

The problem is that this statement says nothing about the
behaviour of the (semi-)stable loci under restriction: we 
are going to investigate this issue.

\begin{proposition}{\label{prop:lin}}
Let $\chi\in\cal X^*(G)$ be a character. There are positive 
integers $c,D_\chi$ depending on $\chi$ such that for any 
$d>D_\chi$ 
we have 
$$
\mbb P(K\oplus V)^\sst\bigl(\Si_{d,\chi^c}\bigr)
=\bbV^\sst(G,\chi)\cup
\raise.5ex\hbox{$\bigl[\bbV^\sst(G,\chi)\sm\pi^{-1}(\hat 0)\bigr]$}
\bigr/ 
\lower.5ex\hbox{$G_m$}.
$$
For writing this union, we view $\mbb P(K\oplus V)=\bbV\cup 
\mbb P(V)$. Moreover, 
$$
\bbV\cap\mbb P(K\oplus V)^\s_{(0)}
\bigl(\Si_{d,\chi^c}\bigr)=\bbV^\s_{(0)}(G,\chi). 
$$ 
In particular, if $K[\bbV]^G=K$, it follows that  
$\mbb P(K\oplus V)^\sst\bigl(\Si_{d,{\chi^c}}\bigr)
=\bbV^\sst(G,\chi)$. 
\end{proposition}

\begin{proof}
We choose $c>0$ such that the condition (i) in lemma \ref{lm:grow} 
is fulfilled. Since the right-hand-side of the equality is invariant 
under the change $\chi\leadsto\chi^c$, we may assume that 
$K[\bbV]^{G,\chi}$ is generated by $K[\bbV]^G_\chi$ over $K[\bbV]^G$. 
We choose further $D,D_\chi>0$ as in $\rm (ii_a)$, resp. $\rm (ii_b)$, 
of {\it loc. cit.} 

For $x\in\bbV^\sst(G,\chi)$, there is a homogeneous 
$f\in K[\bbV]^G_\chi$, with $\deg f\leq D_\chi$, such that 
$f(x)\neq 0$. The section 
\begin{align}{\label{eqn:s}}
s\in
\Gamma\bigl(\mbb P(K\oplus V),\eO_{\mbb P(K\oplus V)}(d)\bigr)
\text{ defined by }
s[a,x]:=\bigl[(a,x),a^{d-D_\chi}f(x)\bigr]
\end{align}
is clearly $G$-invariant and does not vanish at $[1,x]$. 

Assume that $K[\bbV]^G\neq K$, and consider 
$x\in\bbV^\sst(G,\chi)\sm\pi^{-1}(\hat 0)$; we have proved that 
there is $f\in K[\bbV]^G_{\chi^D,(dD)}$, with $f(x)\neq 0$. The 
section 
$$
s\in
\Gamma\bigl(\mbb P(K\oplus V),\eO_{\mbb P(K\oplus V)}(dD)\bigr)
\text{ defined by }
s[a,x]:=\bigl[(a,x),f(x)\bigr]
$$
is again $G$-invariant and does not vanish at $[0,x]$. 

For the converse, let us notice that an immediate consequence of 
the lemma \ref{lm:lin} is that 
$\bbV\cap\mbb P(K\oplus V)^\sst(\Si_{d,\chi})
\subset\bbV^\sst(G,\chi)$. 
Also, in the case where $K[\bbV]^G\neq K$, the converse to the 
condition $\rm (ii_a)$ in lemma \ref{lm:grow} says precisely that 
$$
\mbb P(V)\cap\mbb P(K\oplus V)^\sst\bigl(\Si_{d,\chi}\bigr)
\subset 
\bigl[\bbV^\sst(G,\chi)\sm\pi^{-1}(\hat 0)\bigr]\bigr/ G_m
$$

In the case $K[\bbV]^G=K$, we must prove that the hyperplane 
$\mbb P(V)\hra \mbb P(K\oplus V)$ consists of $\Si_{d,\chi}$-unstable 
points. Let us assume that $[0,x_0]$ is semi-stable; then there is 
$N>0$ and $F\in\Sym^{Nd}\bigl((K\oplus V)^\vee\bigr)$ such that 
$$
F(a,\Si(g,x))=\chi^{N}(g)F(a,x)\text{ and }F(0,x_0)\neq 0. 
$$
Consider $f\in\Sym^\bullet(V^\vee)$, $f(x):=F(0,x)$. Then 
$F(0,x_0)\neq 0$ implies $\deg f=\deg F=Nd$. On the other hand, 
$f(\Si(g,x))=\chi^{N}(g)f(x)$, that is $f\in K[\bbV]^G_{\chi^N}$. 
Applying $\rm (ii_b)$ of lemma \ref{lm:grow}, we find that 
$\deg f\leq ND_\chi<Nd$, a contradiction. 

For the second equality, notice that the left-hand-side is 
obviously included in the right-hand-side. Conversely, we observe 
that for $x\in\bbV^\s_{(0)}(G,\chi)$, the section defined by 
\eqref{eqn:s} vanishes on the hyperplane 
$\mbb P(V)\hra\mbb P(K\oplus V)$, and therefore the orbit $G[1,x]$ 
is closed in $\mbb P(K\oplus V)^\sst(\Si_{d,\chi})$. 
\end{proof}

\begin{remark}{\label{rmk:1}}
We indicate a procedure which allows, possibly after taking finite 
covers of $G$, to linearize the $G$ action in 
$\euf O_{\mbb P(K\oplus V)}(1)$ without altering the semi-stable 
locus. For $d>0$, we consider the endomorphism 
$\zeta_d:Z(G)^\circ\rar Z(G)^\circ$ which raises the elements to 
the $d^{\rm th}$ power, and define
$$
\xymatrix@C=1.8cm{
m_d:Z(G)^\circ\times [G,G]
\ar[r]^-{\phantom{M}(\zeta_d,{\rm id_{[G,G]}})\phantom{M}}
&
Z(G)^\circ\times [G,G]
\xrar{\phantom{M}m_1\phantom{M}}
G,
}
$$
where $m_1$ is the multiplication. Then 
$m_d^*:\cal X^*(G)\rar\cal X^*\bigl(Z(G)^\circ\times [G,G]\bigr)$ 
has finite cokernel, and for all $\chi\in\cal X^*(G)$ we have 
$\mbb V^\sst(G,\chi)=\mbb V^\sst(Z(G)^\circ\times[G,G],m_d^*\chi)$. 
Moreover, the representation 
$$
Z(G)^\circ\times[G,G]\rar\Gl(K\oplus V),\quad
(z,g')
\lmt
{\rm diag}\bigl(\chi^{-1}(z),\chi^{-1}(z)\cdot
(\rho\circ m_d)(z,g')\bigr)
$$ 
induces a linearization of the $Z(G)^\circ\times[G,G]$-action 
$\bar\Si\circ m_d$ in $\crl O_{\mbb P(K\oplus V)}(1)$. 
\end{remark}

In the next section we will compare the numerical criteria for 
semi-stability in the affine and in the projective case, in order 
to address the issue of the chamber structure on $\cal X^*(G)$ 
with respect to the GIT-equivalence relation. 

We fix a norm on $\cal X_*(T)_{\mbb R}$ invariant under the action 
of the Weyl group $W$; it naturally induces a norm $|\cdot|$ on 
$\cal X_*(G)$ (see \cite[subsection 1.1.3]{dh}). 
Consider $x\in\bbV$ and $\l\in\cal X_*(G)$: an immediate computation 
shows that   
$$
\bar\Si\bigl(\l(t),[1,x]\bigr)
=\bigl[1,(t^jx_j)_j\bigr]
\xrar{t\rar 0}
\left\{
\begin{array}[c]{ll} 
[0,x_{m(v,\l)}]&\text{if }m(x,\l)<0;
\\[1.5ex] 
[1,x_0]&\text{if }m(x,\l)=0;
\\[1.5ex]  
[1,0]&\text{if }m(x,\l)>0.
\end{array}
\right.
$$
The weight of the $\l$-action at the specialization, corresponding 
to $\Si_{d,\chi}$, is 
$$
\mu_{d,\chi}(x,\l)=
\left\{
\begin{array}[c]{ll}
\lan\chi,\l\ran-d\cdot m(x,\l)
&
\text{if }m(x,\l)<0;\\[1.5ex] 
\lan\chi,\l\ran
&
\text{if }m(x,\l)\geq 0.
\end{array}
\right.
$$
We consider the obvious extension of $\mu_{\cdot,\cdot}(x,\l)$ to 
$\mbb R\oplus\cal X^*(G)_{\mbb R}$ and define the numerical functions 
\begin{align}{\label{eq:M}}
\begin{array}{l}
\disp
M:\cal X^*(G)_{\mbb R}\times\bbV\rar\mbb R,\quad 
M(\chi,x):=\inf\biggl\{
\frac{\lan\chi,\l\ran}{|\l|}
\,\biggl|\,
\lower0.7ex\hbox{$
\begin{array}{l}
\l\in\cal X_*(G)\\[-1.5ex] 
m(x,\l)\geq 0\biggr.
\end{array}
$}
\biggr\},
\\[1.7ex] 
\disp
\bar M:\bigl(\mbb R\oplus\cal X^*(G)_{\mbb R}\bigr)\times\bbV\rar\mbb R,
\quad
\bar M\bigl((d,\chi),x\bigr):=
\inf_{\l\in\cal X_*(G)}\frac{\mu_{d,\chi}(x,\l)}{|\l|}.
\end{array}
\end{align}
If we could control the locus where $\bar M$ is positive in terms 
of the corresponding locus of $M$ (the converse is easy, since 
$\bar M((d,\cdot),\cdot)\leq M(\cdot,\cdot)$ for $d>0$), the fan 
structure on $\cal X^*(G)_{\mbb R}$ was immediate. Unfortunately, 
at this point we are able to do this only for rational coefficients. 
The equality of the positive loci of $M$ and $\bar M$ follows in 
\ref{cor:numer} as a corollary to our discussion in the next section. 

We conclude with a generality about $G$-actions on 
affine varieties. For a normal, affine, $G$-variety $X$, 
$\Pic^G(X)\cong\cal X^*(G)$, where the isomorphism is given by 
assigning to the character $\chi$ the linearization $\Si_\chi$ 
defined as in \eqref{si-chi}. The following proposition reduces 
the issue of computing the semi-stability subsets of $X$ 
corresponding to $\{\Si_\chi\}_{\chi\in\cal X^*(G)}$ to the 
similar issue for an affine space. 

\begin{proposition}{\label{prop:affine}}
Let $X$ be a normal, affine variety, and $\Si:G\times X\rar X$ a 
regular action on it. Then there is a $G$-module $V$ and a 
$G$-equivariant embedding $\jmath_X:X\hra\bbV:=\Spec(\Sym^\bullet V)$ 
with the property that for all $\chi\in\cal X^*(G)$ holds 
$$
X^\sst(G,\chi)=X\cap\bbV^\sst(G,\chi).
$$
\end{proposition}

\begin{proof}
Let $G':=[G,G]$ be the commutator subgroup of $G$: it is a normal, 
semi-simple subgroup of $G$, and $G/G'$ is a torus. We remark that for 
any character $\chi\in\cal X^*(G)$, $K[X]^G_\chi\subset K[X]^{G'}$ 
and moreover 
$$
\bigoplus_{\chi\in\cal X^*(G)}\kern-1ex 
K[X]^G_\chi
=
\bigoplus_{\hat\chi\in\cal X^*(G/G')}\kern-1ex 
\bigl(K[X]^{G'}\bigr)^{G/G'}_{\hat\chi}
=
K[X]^{G'}.
$$
Since $K[X]^{G'}$ is a finitely generated $K$-algebra and also 
a locally finite $G$-module, there is a finite set 
$\cal S_0\subset K[X]^{G'}$ with the following two properties: 

-- $\cal S_0$ generates $K[X]^{G'}$ as a $K$-algebra; 

-- for each $f\in\cal S_0$ there is $\chi_f\in\cal X^*(G)$ such 
that $f\in K[X]^G_{\chi_f}$. 

\nit We deduce that $V_0:=\sum_{f\in\cal S_0}Kf$ is a $G$-module, 
the natural ring homomorphism $\Sym^\bullet V_0\rar K[X]^{G'}$ is 
an epimorphism of $G$-modules, and moreover 
$(\Sym^\bullet V_0)^G_\chi\rar K[X]^G_\chi$
is surjective for all $\chi\in\cal X^*(G)$. 

Since $K[X]$ is a locally finite $G$-module, we find a finite 
set $\cal S\subset M$ containing $\cal S_0$ such that 

-- $\cal S$ generates $K[X]$ as a $K$-algebra; 

-- $V:=\sum_{f\in\cal S}Kf$ is a $G$-module. 

\nit The natural ring homomorphism $\Sym^\bullet V\rar K[X]$ is 
an epimorphism of $G$-modules, and moreover 
$(\Sym^\bullet V)^G_\chi\rar K[X]^G_\chi$ is surjective for all 
$\chi\in\cal X^*(G)$. A straightforward computation shows that 
the $G$-module $V$ defined in this way fulfills the requirements 
of the proposition. 
\end{proof}


\section{The chamber structure on the cone of effective characters}
{\label{sct:fan}}

Two natural questions naturally arise in the study of the quotients 
$\bbV\invq_\chi G$, as $\chi\in\cal X^*(G)$ varies:

(i) Is there a chamber structure on $\cal X^*(G)$ corresponding to 
the possible semi-stable loci $\bbV^\sst(G,\chi)$?

(ii) For which characters $\chi$ is the corresponding semi-stable 
locus non-empty?\smallskip

\nit When $G$ is a torus, the answer to these questions is given 
by the so-called Gelfand-Kapranov-Zelevinskij (or secondary) fan.

For answering our first question, we will rely 
on the results obtained in \cite{re}, where the author studies the 
GIT-equivalence relation on the space of $G$-linearized line bundles 
over a {\em projective} variety. Let us recall the main result of 
that paper: consider a normal, projective variety $X$, acted on by a 
reductive, linear algebraic group $G$. We define $\NS^G(X)$ to be 
the quotient of $\Pic^G(X)$ by the $G$-algebraic equivalence relation 
(see \cite[subsection 2.1]{re}). It is a finitely generated module 
over $\mbb Z$, so that $\NS^G(X)_{\mbb R}$ is a finite dimensional 
vector space over $\mbb R$. 

One chooses a $W$-invariant norm on $\cal X_*(T)$, and using it 
defines the function $\bar M:\NS^G(X)\times X\rar\mbb R$ such that 
for all $l\in\NS^G(X)$ holds $X^\sst(l)=\{\bar M(l,\cdot)\geq 0\}$. 
Further, one proves that $\bar M$ uniquely 
extends to a function on $\NS^G(X)_{\mbb R}\times X$, which is 
continuous in the first argument (see \cite[lemma 3.2.5]{dh}). One 
extends the definition of the $l$-semi-stable locus to every 
$l\in\NS^G(X)_{\mbb R}$ by the formula $X^\sst(l):=
\{\bar M(l,\cdot)\geq 0\}$. 

\begin{definition}{\label{def:eff}}
(see \cite[section 2]{re}) 
An element $l\in\NS^G(X)_{\mbb R}$ is called {\em effective} if 
$X^\sst(l)\neq\emptyset$. One says that 
$l_1,l_2\in\NS^G(X)_{\mbb R}$ are GIT-{\it equivalent} if 
$X^\sst(l_1)=X^\sst(l_2)$.
\end{definition}
 
The main result of \cite{re} reads:\smallskip

\nit{\sl Theorem}\hskip2ex{\it 
Let $X$ be a normal, projective $G$-variety, and denote by 
$\cal C^G(X)$ the set of effective classes in $\NS^G(X)_{\mbb R}$. 
Then the following hold: 

\nit{\rm (i)} $\cal C^G(X)$ is closed in $\NS^G(X)_{\mbb R}$.

\nit{\rm (ii)} For all $l_o\in\cal C^G(X)$, 
$C(l_o):=\{l\in\cal C^G(X)\mid X^\sst(l_o)\subseteq X^\sst(l)\}$ 
is a closed, convex, \unbar{rational} polyhedral cone in 
$\NS^G(X)_{\mbb R}$. 

\nit{\rm (iii)} The cones $C(l)$, $l\in\cal C^G(X)$ form a fan 
covering $\cal C^G(X)$.  

\nit{\rm (iv)} The {\rm GIT}-equivalence classes are the relative 
interiors of these cones.
}\smallskip

\nit The fan constructed in this way is called in \cite{re} the 
{\rm GIT}-fan, and will be denoted by $\Delta^G(X)$. Further, we 
let $\Delta^G(X)_{\mbb Q}:=\Delta^G(X)\cap\NS^G(X)_{\mbb Q}$ and 
notice that as the cones of $\Delta^G(X)$ are rational, 
$\Delta^G(X)$ is obtained from the fan $\Delta^G(X)_{\mbb Q}$ 
by extending the coefficients from $\mbb Q$ to $\mbb R$. 

What we are going to prove is that this result, which holds for 
projective varieties, can be adapted to our setting, where we 
deal with group actions on affine spaces. For $\chi\in\cal X^*(G)$, 
we define $\bbV^\sst(G,\chi):=\{M(\chi,\cdot)\geq 0\}$, where $M$ 
is given by \eqref{eq:M}. 

\begin{theorem}{\label{thm:fan}}
Let $\rho:G\rar\Gl(V)$ be a representation of the reductive group 
$G$ which has finite kernel, and consider the induced $G$-action on 
$\mbb P(K\oplus V)$ defined by \eqref{eqn:action}. We denote by 
$\re$ the trivial character of $G$. Then the following 
statements hold: 

\nit{\rm (i)} $\mbb R_{\geq 0}(1,\re)$ is a ray in 
$\Delta^G\bigl(\mbb P(K\oplus V)\bigr)$. 

\nit{\rm (ii)} The {\rm GIT}-equivalence classes in 
$\cal X^*(G)_{\mbb R}$ corresponding to the $G$-action on $\bbV$ 
are the relative interiors of the cones of the fan
$$
\Delta^G(\bbV):=
{\rm star}\bigl(\mbb R_{\geq 0}(1,\re)\bigr)\bigl/\mbb R(1,\re). 
$$
\end{theorem}

This result, together with proposition \ref{prop:affine}, 
immediately yield

\begin{theorem}{\label{prop:fan}}
Let $X$ be a normal, affine $G$-variety. The {\rm GIT}-equivalence 
classes in $\cal X^*(G)_{\mbb R}$ corresponding to the $G$-action 
on $X$ are the relative interiors of the cones of a rational, 
polyhedral fan $\Delta^G(X)$, which we call the {\rm GIT}-fan of $X$. 
\end{theorem}

We will prove theorem \ref{thm:fan} in two steps: first we prove 
the existence of the fan structure on the set of effective characters 
in $\cal X^*(G)_{\mbb Q}$; this part uses the results obtained in the 
previous section. Secondly, we prove that the induced fan in 
$\cal X^*(G)_{\mbb R}$, obtained by extending the coefficients from 
$\mbb Q$ to $\mbb R$, parameterizes the GIT-equivalence classes in 
$\cal X^*(G)_{\mbb R}$. 

\begin{proof}(theorem \ref{thm:fan} with $\mbb Q$ coefficients) 
(i) A direct computation yields 
$\mbb P(K\oplus V)^\sst\bigl(\Si_{1,\re}\bigr)
=\bbV\cup[\bbV\sm\pi^{-1}(\hat 0)]\bigr/G_m$, so that $(1,\re)$ 
is an effective class in $\NS^G\bigl(\mbb P(K\oplus V)\bigr)=
\Pic^G\bigl(\mbb P(K\oplus V)\bigr)\cong \mbb Z\oplus \cal X^*(G)$. 
From the third part of the theorem above, we deduce that $(1,\re)$ 
is in the relative interior of a cone in 
$\Delta^G\bigl(\mbb P(K\oplus V)\bigr)_{\mbb Q}$. 

We claim that $\mbb Q_{\geq 0}(1,\re)$ is actually a ray in the 
GIT-fan. Assume the contrary, that $(1,\re)$ is in the relative 
interior of a cone 
$\tau_0\in\Delta^G\bigl(\mbb P(K\oplus V)\bigr)_{\mbb Q}$; 
since $\tau_0$ is rational, we find an integral point $l=(d,\chi)$ 
on its boundary. Then \cite[proposition 8]{re} implies 
that 
$$
\bbV\subset\mbb P(K\oplus V)^\sst\bigl(\Si_{1,\re}\bigr)
\varsubsetneq 
\mbb P(K\oplus V)^\sst\bigl(\Si_{d,\chi}\bigr). 
$$
In particular, 
$[1,0]\in\mbb P(K\oplus V)^\sst\bigl(\Si_{d,\chi}\bigr)$, 
so that there is a homogeneous polynomial $F$ of positive 
degree, such that 
$$
F(1,0)\neq 0\text{ and }
F(a,\rho(g)v)=\chi^N(g)\cdot F(a,v),
\text{ for all }(a,v)\in \bbA^1_K\times\bbV. 
$$
Restricting to $a=1$, we obtain $f\in K[\bbV]^G_{\chi^N}$ such 
that $f(0)\neq 0$. We deduce that    
$$
f(0)=f(\rho(g)0)=\chi^N(g)\cdot f(0),
$$
and therefore $\chi=\re$. It follows that $l=(d,\re)$ is colinear 
with $(1,\re)$, and this contradicts the choice of $l$. 

(ii) First, consider $\chi\in\cal X^*(G)$ such that 
$\bbV^\sst(G,\chi)\neq\emptyset$. For $c,d$ appropriate 
$$
\mbb P(K\oplus V)^\sst\bigl(\Si_{d,\chi^c}\bigr)
\srel{\ref{prop:lin}}{=} 
\bbV^\sst(G,\chi)
\cup
\bigl[\bbV^\sst(G,\chi)\sm\pi^{-1}(\hat 0)\bigr]\bigr/ G_m
\subseteq \mbb P(K\oplus V)^\sst\bigl(\Si_{1,\re}\bigr), 
$$
so that, on one hand, $(d,\chi^c)$ is effective and therefore 
contained in the relative interior of a cone $\tau$ of 
$\Delta^G\bigl(\mbb P(K\oplus V)\bigr)$; on the other hand, 
\cite[proposition 8]{re} implies that $(1,\re)$ belongs to $\tau$. 
All together shows that $\tau\in{\rm star}
\bigl(\mbb Q_{\geq 0}(1,\re)\bigr)$. 

Using proposition \ref{prop:lin} again, we deduce that for two 
GIT-equivalent characters $\chi_1,\chi_2\in \cal X^*(G)$, that 
is $\bbV^\sst(G,\chi_1)=\bbV^\sst(G,\chi_2)$, the linearizations 
$\Si_{d,\chi_1^c}$ and $\Si_{d,\chi_2^c}$ are still 
GIT-equivalent, for suitable $c$ and $d$ sufficiently large. 
Therefore $(d,\chi_1^c)$ and $(d,\chi_2^c)$ are in the relative 
interior of the same cone of the GIT-fan of $\mbb P(K\oplus V)$. 

Conversely, consider 
$\tau\in{\rm star}\bigl(\mbb Q_{\geq 0}(1,\re)\bigr)$ and an 
integral point $(d,\chi)$ in its relative interior. Then lemma  
\ref{lm:lin} implies that 
$\emptyset\neq\bbV\cap\mbb P(K\oplus V)^\sst\bigl(\Si_{d,\chi}\bigr)
\subset \bbV^\sst(G,\chi)$, and therefore $\chi$ is effective. 
Further, since $(1,\re)\in\tau$, $(d+n,\chi)$ is in the relative 
interior of $\tau$ for all $n\geq 0$, and therefore the $(d,\chi)$- 
and $(d+n,\chi)$-semi-stable loci coincide. On the other hand, 
$$
\mbb P(K\oplus V)^\sst\bigl(\Si_{d_1,\chi^{c_1}}\bigr)
=\bbV^\sst(G,\chi)
\cup
\bigl[\bbV^\sst(G,\chi)\sm\pi^{-1}(\hat 0)\bigr]\bigr/ G_m
$$
for some appropriate $c_1$, and for all $d_1$ sufficiently large. 
Choosing $d_1$ and $n$ in such a way that $(d+n,\chi)$ 
and $(d_1,\chi^{c_1})$ are colinear, we deduce that 
$\mbb P(K\oplus V)^\sst\bigl(\Si_{d,\chi}\bigr)$ is still given 
by the right-hand-side of the equality above. 

Consider now $\tau\in{\rm star}\bigl(\mbb Q_{\geq 0}(1,0)\bigr)$, 
and two integral points $(d_1,\chi_1),(d_2,\chi_2)$ in its relative 
interior. Then a similar argument shows that 
$\bbV^\sst(G,\chi_1)=\bbV^\sst(G,\chi_2)$, that is $\chi_1$ and 
$\chi_2$ are GIT-equivalent.
\end{proof}

For extending the coefficients from $\mbb Q$ to $\mbb R$, we have 
to adapt some of the intermediate results in \cite{re} to our context. 
We recall that $T\subset G$ is a maximal torus of $G$; for $g\in G$, 
we denote $T^g:=gTg^{-1}$. Further, the $T$-module $V$ decomposes 
into the direct sum of its weight subspaces $V=\oplus_{a\in\Phi}V_a$, 
and we denote $\{\eta_a\}_{a\in\Phi}$ the corresponding weights. 
Then $\bbV=\times_{a\in\Phi}\bbV_a$, and $T$ acts on $\bbV_a$ by the 
character $\eta_a$. For $x\in\bbV$, we write $(x_a)_{a\in\Phi(x)}$ 
its non-zero coordinates with respect to this decomposition. 

\begin{definition}{\label{stb-set}}
We define the {\em stability set} of a point $x\in\bbV$ to be 
$$
\O(x):=\{\chi\in\cal X^*(G)_{\mbb R}\mid M(\chi,x)\geq 0\}.
$$
\end{definition}

\begin{lemma}{\label{lm:stab-set}}
{\rm (i) }For all $x\in\bbV$, the stability set 
$\O(x)\subset\cal X^*(G)_{\mbb R}$ is a closed, convex, rational 
polyhedral cone.

\nit {\rm (ii)} There are only finitely many stability sets. 
\end{lemma}

\begin{proof} 
(i) Let $x\in\bbV$: $\O(x)$ is closed since $M(\cdot,x)$ is 
continuous on $\cal X^*(G)_{\mbb R}$; $\O(x)$ is convex since 
$M(\chi_1+\chi_2,x)\geq M(\chi_1,x)+M(\chi_2,x)$, for all 
$\chi_1,\chi_2\in\cal X^*(G)_{\mbb R}$. 

We prove now the rationality property. Define the set 
$\cal L_x:=\{\l\in\cal X_*(G)\mid m(x,\l)\geq 0\}$, and 
notice that 
$$
\cal L_x=\bigcup_{g\in G}\cal X_*(T^g)\cap\cal L_x
=\bigcup_{g\in G}
\Ad_g\bigl(\cal X_*(T)\cap \cal L_{\Si(g^{-1},x)}\bigr). 
$$
Moreover, for $g\in G$,  
$\cal X_*(T)\cap \cal L_{\Si(g^{-1},x)}\subset\cal X_*(T)$ 
is convex (if we consider $\mbb Q$-coefficients), and 
$$
\begin{array}{ll}
\cal X_*(T)\cap \cal L_{\Si(g^{-1},x)}
&
=
\{ 
\l\in\cal X_*(T)\mid m\bigl(\Si(g^{-1},x),\l\bigr)\geq 0
\}
\\[1.5ex] 
&\disp
=
\bigl\{
\l\in\cal X_*(T)\,\bigl|\, \min_{a\in\Phi(\Si(g^{-1},x))}
\lan\eta_a,\l\ran\geq 0
\bigr\}.
\end{array}
$$
Since $\Phi(\Si(g^{-1},x))\subset\Phi$ for all $g\in G$, only a 
finite number of such sets appear as $g\in G$ varies; let 
$\Gamma_x\subset G$ be a set of representatives. We deduce that 
$$
\begin{array}{ll}
M(\chi,x)
&\disp
=\inf_{\l\in\cal L_x}\frac{\lan\chi,\l\ran}{|\l|}
=\inf_{g\in G}
\inf
\biggl\{\frac{\lan\chi,\l\ran}{|\l|}
\,\biggl|\,
\l\in\Ad_g(\cal X_*(T)\cap \cal L_{\Si(g^{-1},x)})\biggr\}
\\[1.5ex]
&\disp
=\inf_{g\in G}
\inf
\biggl\{\frac{\lan\chi,\l\ran}{|\l|}
\,\biggl|\,
\l\in\cal X_*(T)\cap \cal L_{\Si(g^{-1},x)}\biggr\}
\\[1.5ex]
&\disp
=\min_{g\in \Gamma_x}
\inf
\biggl\{\frac{\lan\chi,\l\ran}{|\l|}
\,\biggl|\,
\l\in\cal X_*(T)\cap \cal L_{\Si(g^{-1},x)}\biggr\}
\end{array}
$$
and therefore
$$
\begin{array}{ll}
\chi\in\O(x)
&\disp
\Longleftrightarrow
\inf
\biggl\{\frac{\lan\chi,\l\ran}{|\l|}
\,\biggl|\,
\l\in\cal X_*(T)\cap \cal L_{\Si(g^{-1},x)}\biggr\}\geq 0,
\ \forall g\in\Gamma_x
\\[2ex]
&\disp
\Longleftrightarrow
\chi\in\bigcap_{g\in\Gamma_x}
\bigl(
\cal X_*(T)\cap \cal L_{\Si(g^{-1},x)}
\bigr)^\vee.
\end{array}
$$
This is a finite intersection of half-spaces defined by rational 
equations, so that $\O(x)$ is a rational, polyhedral cone. 

(ii) Clearly, any stability set is the union of the GIT classes 
contained in it. We claim that there is a finite number of stability 
sets of the form $\O(x)\cap\cal X^*(G)_{\mbb Q}$: indeed, it follows 
from the proof of theorem \ref{thm:fan} with $\mbb Q$ coefficients 
that 
$$
\O(x)\cap\cal X^*(G)_{\mbb Q}
=
\pr_{\cal X^*(G)_{\mbb Q}}
\bigl[
\O([1,x])\cap 
\text{Support of }{\rm star}\bigl(\mbb Q_{\geq 0}(1,\re)\bigr)
\bigr],
$$
which is clearly finite. By the previous step, 
$\O(x)=\ovl{\O(x)\cap\cal X^*(G)_{\mbb Q}}$, which concludes 
the lemma.  
\end{proof}

\begin{lemma}{\label{chi-cont}}
{\rm (i)} For a subset $U\subset\bbV$, we define 
$C(U):=\{l\in\cal X^*(G)_{\mbb R}\mid U\subseteq
\bbV^\sst(G,l)\}$. Then $C(U)$ is a closed, convex, rational 
cone in $\cal X^*(G)_{\mbb R}$. 

\nit{\rm (ii)} For any $\chi_0\in\cal X^*(G)_{\mbb Q}$ holds 
$$
\bbV^\sst(G,\chi)=\bbV^\sst(G,\chi_0),\quad
\forall 
\chi\in{\rm rel.int.\,}C\bigl(\bbV^\sst(G,\chi_0)\bigr).
$$
\end{lemma}

\begin{proof}
(i) We notice that $C(U)=\bigcap_{x\in U}\O(x)$; by the 
previous lemma, this is a finite intersection of closed, convex, 
rational cones, so that $C(U)$ is the same.

(ii) As $\chi\in C\bigl(\bbV^\sst(G,\chi_0)\bigr)$, it follows 
$\bbV^\sst(G,\chi_0)\subset\bbV^\sst(G,\chi)$, and therefore  
$C\bigl(\bbV^\sst(G,\chi)\bigr)\subset
C\bigl(\bbV^\sst(G,\chi_0)\bigr)$. We deduce that  
$$
\chi\in C\bigl(\bbV^\sst(G,\chi)\bigr)
\cap{\rm rel.int.\,}C\bigl(\bbV^\sst(G,\chi_0)\bigr),
$$
and, using the first part of the lemma, we conclude that there 
is a {\em rational} point $\chi'$ in the intersection above. 
Then $\bbV^\sst(G,\chi)\subset\bbV^\sst(G,\chi')$ since 
$\chi'\in C\bigl(\bbV^\sst(G,\chi)\bigr)$, and also 
$C\bigl(\bbV^\sst(G,\chi')\bigr)=
C\bigl(\bbV^\sst(G,\chi_0)\bigr)$ by theorem \ref{thm:fan} with 
rational coefficients. It follows that 
$\bbV^\sst(G,\chi)\subset\bbV^\sst(G,\chi_0)$. 
\end{proof}

\begin{proof}(theorem \ref{thm:fan} with $\mbb R$ coefficients)
We have already proved the statement for rational coefficients, 
and we know from \cite{re} that the GIT classes in 
$\Delta^G\bigl(\mbb P(K\oplus V)\bigr)$ build a fan structure. 
Combining these facts with lemma \ref{chi-cont}(ii), we deduce 
that for concluding the theorem we still need to prove the following 
statement: any $\chi\in\cal X^*(G)_{\mbb R}$ is GIT-equivalent 
to some rational character $\chi_0\in\cal X^*(G)_{\mbb Q}$. 

Let $\chi\in\cal X^*(G)_{\mbb R}$ such that $\bbV^\sst(G,\chi)\neq 
\emptyset$. We know from lemma \ref{chi-cont} that 
$C\bigl(\bbV^\sst(G,\chi)\bigr)\subset\cal X^*(G)_{\mbb R}$ is 
a rational polyhedral cone. Therefore we find a sequence 
${\{\chi_n\}}_n\subset\cal X^*(G)_{\mbb Q}$ such that 
$\chi_n\rar\chi$. Since the number of stability sets for rational 
classes is finite, we may assume, after possibly passing to a 
subsequence, that $\bbV^\sst(G,\chi_n)=:U$ is independent of $n$. 
Now, $\chi_n\in C\bigl(\bbV^\sst(G,\chi_n)\bigr)=C(U)$ for 
all $n$; since $\chi_n$ converges to $\chi$ and $C(U)$ is closed, 
we deduce that $\chi\in C(U)$. 

There is a unique face $\tau\subset C(U)$, possibly equal to $C(U)$, 
such that, on one hand, $\chi\in{\rm rel.int.\,}\tau$; on the other 
hand, applying again theorem \ref{thm:fan} with rational coefficients, 
we deduce that such a face equals  $C\bigl(\bbV^\sst(G,\chi_0)\bigr)$, 
for any rational point $\chi_0\in{\rm rel.int.\,}\tau$. Applying 
lemma \ref{chi-cont}, we deduce that $\chi$ is GIT-equivalent to 
$\chi_0$. 
\end{proof}

An immediate consequence of the theorem is the following comparison 
result between the numerical functions defined by \eqref{eq:M}:

\begin{corollary}{\label{cor:numer}}
For any $\chi\in\cal X^*(G)$, there is a positive number $d_\chi>0$ 
depending on $\chi$ such that for any $d>d_\chi$ holds: 
$$
M(\chi,x)\geq 0\;(\text{\rm resp.}>0)\Longleftrightarrow 
\bar M\bigl((d,\chi),[1,x])\geq 0\;(\text{\rm resp.}>0),
\quad\forall x\in\bbV.
$$
\end{corollary} 

\nit In this statement, the non-obvious implication is from left to 
the right. Our sign differs from that in \cite{dh} because there the 
authors identify a vector space with the affine space defined by it 
(see \cite[subsection 1.1.5]{dh}). 

Now we turn to our second question, namely to describe those characters 
of $G$ for which the semi-stable locus is non-empty. Is immediate to 
see that the effective characters form a convex cone in 
$\cal X^*(G)_{\mbb R}$, and the natural question which raises is 
what this cone looks like. Secondly, we wish to characterise the 
characters $\chi\in\cal X^*(G)_{\mbb R}$ for which 
$\bbV^\sst(G,\chi)\neq \bbV^\s_{(0)}(G,\chi)$. In the terminology of 
\cite{dh}, one says that these characters {\em belong to a wall}\,; 
on the other hand, in \cite{re} a wall is by definition a cone of 
codimension one in $\Delta^G(V)$. Here we will adopt the former 
definition, but we stress that there may be cones of codimension zero 
in $\Delta^G(\bbV)$ which are walls for this definition. 

\begin{definition}{\label{def:bad}}
A character $l\in\cal X^*(G)_{\mbb R}$ is {\em effective} if 
$\bbV^\sst(G,l)\neq\emptyset$. 
A cone $\tau\in\Delta^G(\bbV)$ is a {\em wall} if $\bbV^\sst(G,l)
\neq \bbV^\s_{(0)}(G,l)$ for all $l\in{\rm rel.int.\,}\tau$. 
\end{definition}

We start by addressing these issues in the abelian case. 

\begin{lemma}{\label{prop:eff}}
Let $G$ be a torus and $\rho:G\rar\Gl(V)$ a representation, whose 
weights are $\{\eta_a\}_{a\in\Phi}$. 

\nit{\rm (i)} The effective cone is 
$\;\disp\cal C^G(\bbV)=\sum_{a\in\Phi}\mbb R_{\geq 0}\eta_a.$

\nit{\rm (ii)} The walls are precisely the cones 
$$
\tau_{_{\Phi'}}:=\sum_{a\in\Phi'}\mbb R_{\geq 0}\eta_a, 
\text{ where } \Phi'\subset\Phi
\text{ is such that }
\codim_{\cal X^*(G)_{\mbb R}}\lan\eta_a;a\in\Phi'\ran=1.
$$
\end{lemma}

\begin{proof} 
(i) Since the GIT-fan is rational, we may restrict ourselves to 
integral coefficients. Consider 
$\chi=\sum_{a\in\Phi(\chi)}k_a\eta_a$, with $\Phi(\chi)\subset\Phi$ 
and $k_a>0$. For each $a\in\Phi(\chi)$, we consider a linear function 
$f_a$ on the weight space $\bbV_a\hra\bbV$. Then 
$f:=\prod_{a\in\Phi(\chi)}f_a^{k_a}\in K[\bbV]$ is $\chi$-invariant, 
and not identically zero. More precisely, $f$ does not vanish 
at general points of the form $x=(x_a)_a$, with $x_a=0$ for 
$a\not\in\Phi(\chi)$. 

Conversely, let $\chi\in\cal X^*(G)$ be effective, and choose  
$x=(x_a)_{a\in\Phi(x)}\in \bbV^\sst(G,\chi)\sm\{0\}$. Using 
proposition \ref{prop:sst}, we deduce that 
\begin{align}{\label{eq:ss}}
\begin{array}{lcl}
x\text{ is $\chi$-semi-stable}
&
\Longleftrightarrow
&
\biggl[
\forall\,\l\in\{\lan\chi,\cdot\ran<0\}\Rightarrow
\min_{a\in\Phi(x)}\{\lan\eta_a,\l\ran\}<0
\biggr]
\\ 
&
\Longleftrightarrow
&
\{\lan\chi,\cdot\ran<0\}\cap
\left(\sum_{a\in\Phi(x)}\mbb Z_{\geq 0}\eta_a\right)^{\!\vee}
=\emptyset
\\[1ex] 
&
\Longleftrightarrow
&
\left(\sum_{a\in\Phi(x)}\mbb Z_{\geq 0}\eta_a\right)^{\!\vee}
\subset
\{\lan\chi,\cdot\ran\geq 0\}
\\[1ex] 
&
\Longleftrightarrow
&
\chi\in\sum_{a\in\Phi(x)}\mbb Z_{\geq 0}\eta_a.
\end{array}
\end{align}
Since $\sum_{a\in\Phi(x)}\mbb Z_{\geq 0}\eta_a\subset 
\sum_{a\in\Phi}\mbb R_{\geq 0}\eta_a$, the first statement follows. 

(ii) A character $\chi\in\cal X^*(G)$ is in a wall precisely if there 
is a point $x=(x_a)_{a\in\Phi(x)}\in\bbV^\sst(G,\chi)$ with positive 
dimensional stabilizer. We have proved in \eqref{eq:ss} that this 
implies $\chi\in\sum_{a\in\Phi(x)}\mbb Z_{\geq 0}\eta_a$. Also, there 
is a 1-PS $\l$ contained in the stabilizer of $x$, which amounts 
requirering  that $\lan\eta_a,\l\ran=0$ for all $a\in\Phi(x)$. Therefore 
the $\{\eta_a\}_{a\in\Phi(x)}$'s does not span $\cal X^*(G)_{\mbb R}$. 
Conversely, one immediately sees that the point constructed in the 
previous part has positive dimensional stabilizer. 
\end{proof}

We consider now the general case of a representation $\rho:G\rar\Gl(V)$
of a reductive group. We denote
$$
\cal X^*(G,\bbV):=
\left\{
\chi\in\cal X^*(G)
\,\left|\, 
\begin{array}{l}
\mbb R_{\geq 0}\chi\text{ is a ray in }
\Delta^T(\bbV)\cap\cal X^*(G)_{\mbb R}
\\ 
\mbb Z_{\geq 0}\chi= \mbb R_{\geq 0}\chi\cap X^*(G)
\\ 
\chi\text{ is }G\text{-effective}
\end{array}
\right.
\right\}.
$$

\begin{proposition}{\label{prop:wall}}
{\rm (i)} Let $G$ be a reductive group and $\rho:G\rar\Gl(V)$ be a 
representation with finite kernel. The cone of effective characters 
for the action of $G$ on the affine space $\bbV$ is convex and equals 
$$
\cal C^G(\bbV)=\sum_{\chi\in{\cal X}^*(G,\bbV)}\mbb R_{\geq 0}\chi.
$$ 
In the case $K[\mbb V]=K$, $\cal C^G(\bbV)$ is strictly convex. 

\nit {\rm (ii)} Let $G$ be a reductive group and $\rho:G\rar\Gl(V)$ be 
a representation with finite kernel, and let $\{\eta_a\}_{a\in\Phi}$ 
be the weights of its maximal torus. Assume that $\cal X^*(G)_{\mbb R}$ 
is not contained in the linear span of any $T$-wall. Then the 
walls in $\Delta^G(\bbV)$ are of the form 
$\cal X^*(G)_{\mbb R}\cap\tau_{_{\Phi'}}$, with $\Phi'\subset\Phi$ 
as above. In particular, in this case, {\em all} walls have 
codimension one in $\Delta^G(\bbV)$. 
\end{proposition}

\begin{proof} (i) For the same reason as before, we prove the 
statements for integral characters. 
The inclusion of the right-hand-side into the left-hand-side 
is obvious. For the other direction, let us consider 
$\chi\in\cal X^*(G)$ such that $\mbb R_{\geq 0}\chi$ is an exterior 
ray of $\cal C^G(\bbV)$. As $\cal C^G(\bbV)\subset{\rm Support}
\bigl(\Delta^T(\bbV)\cap\cal X^*(G)_{\mbb R}\bigr)$, there is a cone 
$\tau\in\Delta^T(\bbV)$ such that $\chi\in{\rm rel.int.}
(\tau\cap\cal X^*(G))$. Then, for any $\chi'\in\tau\cap\cal X^*(G)$ 
holds 
$$
\bbV^\sst(T,\chi)\subset\bbV^\sst(T,\chi')
\srel{\ref{lm:us}}{\Longrightarrow}
\bbV^\us(G,\chi)\supset\bbV^\us(G,\chi'),
$$
so that $\chi'$ is still $G$-effective. Since $\chi$ generates an 
exterior ray in $\cal C^G(\bbV)$, 
$\mbb R_{\geq 0}\chi=\tau\cap\cal X^*(G)_{\mbb R}$; otherwise would 
be contained in the relative interior of a `larger' cone. It follows 
that $\chi\in\cal X^*(G,\bbV)$.

(ii) We notice that the induced representation $\rho_T:T\rar\Gl(V)$ 
still has finite kernel, and therefore the weights 
$\{\eta_a\}_{a\in\Phi}$ span $\cal X^*(T)_{\mbb R}$. The Hilbert-Mumford 
criterion implies that any wall in $\Delta^G(\bbV)$ is contained in 
the intersection of $\cal X^*(G)_{\mbb R}$ with a wall in 
$\Delta^T(\bbV)$. Our assumption, combined with the fact that the walls 
for the $T$-action on $\bbV$ have codimension one in $\cal X^*(T)$, 
imply that for any wall $\tau\in\Delta^T(\bbV)$, the intersection 
$\cal X^*(G)_{\mbb R}\cap\tau\subset\Delta^G(\bbV)$ has codimension 
one. We use now \cite[proposition 3.3.20]{dh} (which, by theorem 
\ref{thm:fan}, is still valid in our affine setting), which says that 
every interior wall in $\Delta^G(\bbV)$ is contained in a wall of 
codimension {\em at most} one, and deduce that all the interior walls 
in $\Delta^G(\bbV)$ have the form $\cal X^*(G)_{\mbb R}\cap\tau$, with 
$\tau\in\Delta^T(\bbV)$ a wall. 

This argument leaves uncovered the exterior cones of $\Delta^G(\bbV)$, 
but these are obviously walls in the sense of \ref{def:bad}, and they 
have (by definition) codimension one. 
\end{proof}

Let us remark that the very same argument as in 
\cite[proposition 3.2.8]{dh} shows that if $l\in\cal C^G(\bbV)$ 
belongs to an exterior cone of $\Delta^G(\bbV)$, then 
$\bbV^\s_{(0)}(G,l)=\emptyset$.


\section{The Chow ring of the quotients}{\label{sct:chow}}

Assume again that $\rho:G\rar\Gl(V)$ is a representation with 
finite kernel. Our goal in this section is to compute the Chow 
ring of the invariant quotients of $\bbV$, corresponding to 
characters of $G$. The main tool for this 
computation is the equivariant intersection theory developed 
in \cite{eg}. For shorthand, we denote by $A_*^G$ the equivariant 
Chow ring of a point. 

First we recall some results from \cite{dem}: the {\em discriminant} 
$\Delta_G\in\End(\Sym^\bullet\cal X^*(G))$ is defined to be the product 
of the linear forms corresponding to the reflections in the Weyl group. 
This definition determines it is uniquely, up to sign. Equivalently, 
for a choice $T\subset B\subset G$ of a Borel subgroup, $\Delta_G$ 
is defined as the product of the positive roots. For $w\in W=N_G(T)/T$, 
we let $\det(w)$ be the determinant of $w$, viewed as an endomorphism 
of $\cal X^*(G)$. The discriminant has the property that for all 
$w\in W$, $w(\Delta_G)=\det(w)\Delta_G$.
Since $A^*_T\cong\Sym^\bullet\cal X^*(G)$, we will 
view $\Delta_G$ as an element in $A^{\dim G/B}_T$. 

We consider the $(A^*_T)^W$-linear endomorphism 
$\disp J_G:=\sum_{w\in W}\det(w)w$ of $A^*_T$; according to 
\cite[lemme 4]{dem}, it has the property 
that $\wp_G:=J_G/\Delta_G$ is still and endomorphism of $A^*_T$, 
which takes values in $(A^*_T)^W$. Since $\wp_G:(A^*_T)_{\mbb Q}
\rar(A^*_T)^W_{\mbb Q}$ is $(A^*_T)^W_{\mbb Q}$-linear, and 
$\wp_G(\Delta_G)=|W|$, it is an epimorphism. 

The geometrical meaning of the homomorphism $\wp_G$ is captured 
in the 

\begin{lemma}{\label{lm:mean}}
The composition $A^*_T\cong A^*_B\srel{\kern-0.3ex\vphi_{\!*}}{\rar}
A^{*-\dim G/B}_G\srel{\kern0.2ex\vphi^*}{\rar}A^{*-\dim G/B}_T$ 
equals $\wp_G$. 
\end{lemma}

\begin{proof}
The push-forward formula implies that the endomorphism 
$h:=\vphi^*\vphi_*$ is $\vphi^*(A^*_G)$-linear. But 
$\vphi^*:A^*_G\rar (A^*_T)^W$ is injective, its cokernel 
is $\mbb Z$-torsion, and therefore $h$ is $(A^*_T)^W$-linear. 
Moreover, $h$ is vanishes in degree strictly less than 
$\dim G/B=\deg\Delta_G$. Applying \cite[proposition 1]{dem}, 
we obtain that 
$$
|W|\cdot\Delta_G\cdot h(u)=h(\Delta_G)\cdot J_G(u)
=|W|\cdot J_G(u), 
$$
and our claim follows. 
\end{proof}

With these preparations, the main result of this section is 

\begin{theorem}{\label{thm:chow}}
Let $\chi\in\cal X^*(G)$ be a character such that 
$\bbV^\sst(G,\chi)=\bbV^\s_{(0)}(G,\chi)$; in particular 
$\bbV\invq_\chi G$ is a geometric quotient. Then the Chow ring  
$$
A_*\bigl(\bbV\invq_\chi G\bigr)_{\mbb Q}\cong 
\biggl.
\bigl(A_*^T\bigr)^W_{\mbb Q}
\biggr/
\wp_G\bigl\lan [\bbE(\l)]_T;\l\in\euf F(\chi)\bigr\ran_{\mbb Q}.
$$
\end{theorem}

We wish to remark that this result represents a considerable 
generalization of \cite[theorem 4.4]{es}, where one assumes 
that $\rho(G)$ contains the homotheties of $\Gl(V)$ (which implies 
that $K[\bbV]^G=K$), and moreover that $G$ acts freely on 
$\bbV^\s_{(0)}(G,\chi)$. The authors of \cite{eg} were aware of 
the possibility of generalizing the results of Ellingsrud and 
Str\o mme, as they explicitely point this out on page 610 of 
{\it loc.\,cit.} The computation of $[\bbE(\l)]_T$ is immediate: 
if $E\hra V$ is a sub-$T$-module, then the equivariant class of 
the linear space $\bbE:=\Spec\bigl(\Sym^\bullet(V/E)\bigr)$  is
$$
[\bbE]_T=\prod_a\eta_{a}^{m_{a}}\in A_*^T,
$$
where $\bigl\{\eta_{a}\bigr\}_a$ are the weights of the $T$-module 
$E$, and $\bigl\{m_{a}\bigr\}_a$ are the corresponding multiplicities. 

For the proof of this theorem we need some preparatory results. 

\begin{lemma}{\label{lm:surj}}
Let $\si:X\rar Y$ a projective, surjective morphism between 
reduced and irreducible, quasi-projective  varieties. Then 
$\si_*:A_*(X)_{\mbb Q}\rar A_*(Y)_{\mbb Q}$ is surjective. 
\end{lemma}

\begin{proof}
We consider an effective cycle $Z\hra Y$ and prove that there is 
an effective cycle $W\hra X$ such that $\si|_W:W\rar X$ is finite. 
Since $\si$ is surjective, there is a component $X'\hra\si^{-1}(Z)$ 
which maps onto $Z$; the restriction $\si|_{X'}:X'\rar Z$ is still 
projective. Consider an irreducible component $W\hra X'$ of 
a general hyperplane section of $X'$, of codimension equal the 
dimension of the general fibre of $\si|_W$; it will have the 
property that the map $\si|_W:W\rar Z$ is finite. It follows 
that $\si_*[W]$ is a positive multiple of $[Z]$. 
\end{proof}

\begin{lemma}{\label{lm:image}} Let $P\subset G$ be a 
parabolic subgroup, $L\subset P$ its Levi component, and 
consider $\tilde\eta\in A^*_T(G/P)^W$. Further, denote by 
$f_*^T:A^*_T(G/P)\rar A^*_T$ the proper push-forward.  

\nit{\rm (i)} The fixed point set of $G/P$ under the natural 
$T$-action consists of finitely many points; these are 
$$
{\bigl(G/P\bigr)}^T
=\bigl\{
wP\mid w\in N_G(T)/N_L(T)=W/W_L
\bigr\},
$$
\nit{\rm (ii)} Let $\imath_1:\{P\}\hra G/P$ be the inclusion, 
and $\eta:=\imath_1^*\tilde\eta$. Then
$$
f_*^T\bigl(\lan\tilde\eta\ran^W_{\mbb Q}\bigr)=
\wp_G\lan\eta\ran_{\mbb Q}.
$$ 
\end{lemma}

\begin{proof} The first statement is well-known. For the second 
one, we will use the integration formula \cite[corollary 1]{eg1} 
applied to the homogeneous variety $G/P$. For $w\in W$ we have 
the diagram 
$$
\xymatrix@R=0.7cm{
\{wP\}\;
\ar@{^(->}[r]^-{\imath_w}\ar[d]_-{f_w}
&
G/P
\ar[d]^-{f}
\\ 
\Spec K
\ar@{=}[r]
&
\Spec K
}
$$
Observe that for $\tld\alpha\in\lan\tld\eta\ran^W_{\mbb Q}$, 
$\imath_w^*\tld\alpha=w(\imath_1^*\tld\alpha)$; in particular 
$\imath_1^*\tld\alpha\in A^*_T$ is $W_L$-invariant. Choosing a 
set of representatives $\hat W\subset W$ for the rest classes 
in $W/W_L$, with $1\in\hat W$, the integration formula reads 
\begin{align}
f_*^T(\tld\alpha)
=\sum_{\hat w\in\hat W} {(f_{\hat w}^T)}_*
\frac{\imath_{\hat w}^*\tld\alpha}
{\eul^T\bigl({\sf T}_{\hat wP}(G/P)\bigr)}
=\sum_{\hat w\in\hat W} \hat w\biggl(
\frac{{(f_{1}^T)}_*\imath_{1}^*\tld\alpha}
{{(f_{1}^T)}_*\eul^T\bigl({\sf T}_{P}(G/P)\bigr)}
\biggr).
\end{align}
We observe now that the image of ${(f_{1}^T)}_*\imath_{1}^*: 
A^*_T(G/P)^W_{\mbb Q}\rar A^*_T$ equals $(A^*_T)^{W_L}_{\mbb Q}$, 
and also that ${(f_{1}^T)}_*\imath_{1}^*\tld\eta=\eta$. 
It follows that 
${(f_{1}^T)}_*\imath_{1}^*\lan\eul^T(\crl Q)\ran^W_{\mbb Q}
=\eta\cdot (A^*_T)^{W_L}_{\mbb Q}$. 

Secondly, we claim that 
${(f_{1}^T)}_*\eul^T\bigl({\sf T}_{P}(G/P)\bigr)$ equals the quotient 
$\Delta_G/\Delta_L$ of the discriminants of $G$ and $L$ respectively. 
Indeed, if $\mfrak g$, $\mfrak p$ and $\mfrak b$ are the Lie algebras 
of $G,P$ and $B$ respectively, then 
${(f_{1}^T)}_*\eul^T\bigl({\sf T}_{P}(G/P)\bigr)$ is the product of 
weights of the $T$-module $\mfrak g/\mfrak p\cong(\mfrak g/\mfrak b)
\bigr/(\mfrak p/\mfrak b)$. Remains to notice that the weights of 
$T$-module $\mfrak p/\mfrak b$ are (up to sign) precisely the positive 
roots of $L$. 

We rewrite equality above as
\begin{align}{\label{local}}
\begin{array}{ll}
f_*^T(\tld\alpha)
&\disp
=\sum_{\hat w\in\hat W} \hat w\biggl(
\frac{\alpha}{\Delta_G/\Delta_L}
\biggr)
,\quad\text{with }
\alpha:={(f_{1}^T)}_*\imath_{1}(\tld\alpha)\in
\eta\cdot(A^*_T)^{W_L}_{\mbb Q}
\\[1.5ex]
&\disp
=\frac{1}{|W_L|}\sum_{w\in W}w
\biggl(\frac{\alpha}{\Delta_G/\Delta_L}\biggr) 
=\frac{1}{|W_L|}\wp_G(\Delta_L\cdot\alpha).
\end{array} 
\end{align}
We are going to compute $\wp_G\lan\eta\ran_{\mbb Q}$, and 
compare with this formula. For the same reason as before, 
$\eta\in (A^*_T)^{W_L}$. 
For $\eta b\in\lan\eta\ran_{\mbb Q}$, we obtain 
$$
\begin{array}{ll}
\wp_G(\eta b)
&\disp
=\frac{1}{\Delta_G}\sum_{w\in W}\det(w)\cdot w(\eta b)
=\frac{1}{\Delta_G}\sum_{\hat w\in\hat W}\sum_{w\in\hat w W_L}
\det(w)\cdot w(\eta b)
\\[1,5ex]
&\disp 
=\frac{1}{\Delta_G}\sum_{\hat w\in\hat W}
\hat w(\eta)\cdot\sum_{w\in\hat w W_L}
\det(w)\cdot w(b)
\\[1,5ex]
&\disp
=\frac{1}{\Delta_G}\sum_{\hat w\in\hat W}
\det(\hat w)\cdot\hat w(\eta)\cdot\hat w
\biggl(
\sum_{w\in W_L}\det(w)\cdot w(b)
\biggr)
\\[1,5ex]
&\disp 
=\frac{1}{\Delta_G}\sum_{\hat w\in\hat W}
\det(\hat w)\cdot\hat w(\eta)\cdot\hat w
\bigl(\Delta_L\cdot \wp_L(b)\bigr)
=
\sum_{\hat w\in\hat W}\hat w
\biggl(\frac{\eta\cdot\wp_L(b)}{\Delta_G/\Delta_L}\biggr). 
\end{array}
$$
As $\wp_L:(A^*_T)_{\mbb Q}\rar (A^*_T)^{W_L}_{\mbb Q}$ is 
surjective, $\{\eta\cdot\wp_L(b)\mid b\in (A^*_T)_{\mbb Q}\}
=\eta\cdot (A^*_T)^{W_L}_{\mbb Q}$, which concludes the proof. 
\end{proof}

\begin{proposition}{\label{prop:image}}
Let $\rho:G\rar \Gl(V)$ be a representation with finite kernel, 
and $\bbE\hra\bbV$ a linear subspace such that its stabilizer 
$P:={\rm Stab}_G(\bbE)$ is a parabolic subgroup of $G$. Then 
the followings hold:

\nit{\rm (i)} $G\cdot \bbE$ is closed in $\bbV$; 

\nit{\rm (ii)} the image of the canonical homomorphism 
$$
A_*^G(G\cdot \bbE)_{\mbb Q}\lar A_*^G(\bbV)_{\mbb Q}
\cong {(A_{*-\dim \bbV}^G)}_{\mbb Q}
$$
equals $\wp_G\lan[\bbE]_T\ran_{\mbb Q}$. 
\end{proposition}

\begin{proof}
(i) Let $e:=\dim\bbE$, and consider the action of $G$ on the 
Grassmannian $\Grs(e,\bbV)$ of $e$-dimensional linear subspaces of 
$\bbV$. We denote $\crl T\rar\Grs(e,\bbV)$ the tautological bundle 
of rank $e$, and by $\rO_\bbE\cong G/P$ the $G$-orbit of 
$[\bbE]\in\Grs(e,\bbV)$; since $P\subset G$ is parabolic, $\rO_\bbE$ 
is closed in $\Grs(e,\bbV)$. 

In the following diagram
$$
\xymatrix@R=0.7cm{
\crl T|_{\rO_\bbE}
\ar[d]\ar@{^(->}[r]^-{\jmath}
&
\rO_\bbE\times\bbV
\ar[d]\ar[r]^-{\pr_\bbV}
&
\bbV
\ar[d]
\\ 
\rO_\bbE
\ar@{=}[r]
&
\rO_\bbE
\ar[r]^-{f}
&
\Spec K
}
$$
all the morphisms are $G$-equivariant, and we observe that 
$G\cdot\bbE=\pr_\bbV(\crl T|_{\rO_\bbE})$. Since $\crl T|_{\rO_\bbE}$ 
is closed in $\rO_\bbE\times \bbV$ and $\pr_\bbV$ is proper, 
$G\cdot\bbE$ is closed in $\bbV$. 

(ii) From the commutative diagram 
$$
\xymatrix{
\crl T|_{\rO_\bbE}
\ar@{^(->}[r]^-{\jmath}\ar[d]_-{\si}
&
\rO_\bbE\times\bbV
\ar[d]^-{\pr_\bbV}
\\ 
G\cdot\bbE\ 
\ar@{^(->}[r]
&
\bbV
}
$$
and the surjectivity of $\si_*$ proved in \ref{lm:surj}, we deduce 
that image of $A_*^G(G\cdot \bbE)_{\mbb Q}\lar A_*^G(\bbV)_{\mbb Q}$ 
equals the image of ${{(\pr_\bbV)}}_*^G\circ\jmath_*^G$. 
We consider the commutative diagram of homomorphisms between 
$G$-equivariant Chow groups 
$$
\xymatrix@C=2cm@R=0.7cm{
A_*^G(\crl T|_{\rO_\bbE})
\ar[r]^-{\jmath_*^G}
&
A_*^G(\rO_\bbE\times\bbV)
\ar[r]^-{{{(\pr_\bbV)}}_*^G}
&
A_*^G(\bbV)
\\
A_{*-e}^G(\rO_\bbE)
\ar[u]^-\cong\ar[r]^-{\eul^G(\cal Q|_{\rO_\bbE})\,\cap}
&
A_{*-\dim \bbV}^G(\rO_\bbE)
\ar[u]^-\cong\ar[r]^-{f_*^G}
&
A_{*-\dim \bbV}^G
\ar[u]^-\cong
}
$$
and the excess intersection theorem (see \cite[example 6.3.5]{fu}) 
implies that the lower left homomorphism 
is the cap product with the equivariant Euler class of the universal 
quotient bundle $\crl Q\rar \Grs(e,\bbV)$. By abuse of notation, we will 
still write $\crl Q\rar \rO_\bbE$ for the restriction of $\crl Q$ to 
$\rO_\bbE$. Let us consider an appropriate open subset $U$ in some 
representation space of $G$ (see \cite[definition-proposition 1]{eg}) 
needed for the computation of the equivariant Chow groups. The diagram 
$$
\xymatrix@R=0.7cm{
U\times_T\rO_\bbE
\ar[d]_-{\vphi_{U,\rO_\bbE}}\ar[r]^-{f^T}
&
U/T
\ar[d]^-{\vphi_U}
\\ 
U\times_G\rO_\bbE
\ar[r]^-{f^G}
&
U/G
}
$$
is cartesian, the horizontal morphisms are proper and the vertical 
ones are flat, and therefore the induced diagram 
$$
\xymatrix@C=1.5cm@R=0.7cm{
\bigl\lan\eul^T(\crl Q)\bigr\ran_{\mbb Q}\subset
{A_*^T(\rO_\bbE)}_{\mbb Q}
\ar[r]^-{f_*^T}
&
{(A_*^T)}_{\mbb Q}
\\ 
\bigl\lan\eul^G(\crl Q)\bigr\ran_{\mbb Q}\subset
{A_*^G(\rO_\bbE)}_{\mbb Q}
\ar[r]^-{f_*^G}\ar@{->}[u]^-{\vphi^*}
&
{(A_*^G)}_{\mbb Q}
\ar@{->}[u]_-{\vphi^*}
}
$$
commutes (see \cite[proposition 1.7]{fu}). Now, we know from 
\cite[proposition 6]{eg} that $\vphi^*$ is injective, and its 
image consists of the $W$-invariant elements, that 
is 
$$
\vphi^*\bigl[
f_*^G\bigl(\lan\eul^G(\crl Q)\ran_{\mbb Q}\bigr)
\bigr]
=
f_*^T\bigl(\lan\eul^T(\crl Q)\ran^W_{\mbb Q}\bigr)
\srel{\ref{lm:image}}{=}
\wp_G\bigl(\lan[\bbE]_T\ran^W_{\mbb Q}\bigr).
$$
All together shows that the image of 
$A_*^G(G\cdot\bbE)\rar A_{*-\dim \bbV}^G$ is 
$\wp_G\bigl(\lan[\bbE]_T\ran^W_{\mbb Q}\bigr)$. 
\end{proof}

\begin{proof}(of the theorem)
We have proved in lemma \ref{lm:us} that 
$$
\bbV^\s_{(0)}(G,\chi)/G
=\bigl(\bbV\sm \bbV^\us(G,\chi)\bigr)\bigl/G\bigr.
=\biggl(
\bbV\sm\bigcup_{\l\in\euf F(\chi)}G\cdot \bbE(\l)
\biggr)
\biggl/G\biggr.,
$$
where the $\bbE(\l)$'s are linear subspaces of $\bbV$. Applying 
theorem 3 and proposition 5 of \cite{eg}, we deduce that 
$$
A_{*-\dim G}\bigl(\bbV^\s_{(0)}(G,\chi)/G\bigr)\cong
A_*^G\bigl(\bbV^\s_{(0)}(G,\chi)\bigr),
$$
and that there is an exact sequence 
$$
A_*^G\biggl(\bigcup_{\l\in\euf F(\chi)}G\cdot \bbE(\l)\biggr)
\srel{\jmath_*^G}{\lar}
A_*^G(\bbV)
\lar
A_*^G\bigl(\bbV^\s_{(0)}(G,\chi)\bigr)
\lar 0.
$$
Further, by the very construction of $\euf F(\chi)$, the 
$\{G\cdot \bbE(\l)\}_{\l\in\euf F(\chi)}$ are the irreducible 
components of the union, and therefore 
$$
\jmath_*^G
\biggl[
A_*^G\biggl(\bigcup_{\l\in\euf F(\chi)}G\cdot \bbE(\l)\biggr)
\biggr]
=
\sum_{\l\in\euf F(\chi)}
\jmath_*^G
\bigl[
A_*^G\bigl(G\cdot \bbE(\l)\bigr)
\bigr].
$$
For $\l\in\euf F(\chi)$, $\bbE(\l)$ is stabilized by the parabolic 
subgroup $P(\l)\subset G$, so that ${\rm Stab}_G\bbE(\l)$ is parabolic 
too. It follows now from proposition \ref{prop:image} that 
$$
\jmath_*^G
\bigl[
A_*^G\bigl(G\cdot \bbE(\l)\bigr)
\bigr]_{\mbb Q}
=\wp_G\lan [\bbE(\l)]_T\ran_{\mbb Q}.
$$ 
This finishes the proof of theorem \ref{thm:chow}. 
\end{proof}

\begin{remark}{\label{rmk:more}}
We have been primarily interested in the Chow ring of the quotients, 
but the proof shows more: for {\em any} $\chi\in\cal X^*(G)$, the 
$G$-equivariant Chow ring of the $\chi$-(semi-)stable locus of 
$\bbV$ is given by the formula appearing in theorem \ref{thm:chow}. 
We precise that for computing the equivariant ring of the properly 
stable locus, one has to divide out the (possibly) larger ideal 
$\wp_G\bigl\lan 
[\bbE(\l)]\,;\,\l\in\{\lan\chi,\,\cdot\,\ran\leq 0\}
\bigr\ran_{\mbb Q}$. 

The difference between the semi-stable and the properly stable locus 
is particularly relevant if $G$ is semi-simple, when its character 
group is trivial. The semi-stable locus is then the whole $\bbV$, and 
the corresponding invariant quotient is not geometric. On the other 
hand the quotients like $\bbV^\s_{(0)}(G)/G$ can be used to construct 
approximations of the algebro-geometric classifying space of $G$. 
\end{remark}

We conclude this section with a remark about the generators of the 
Chow ring. In \cite[section 6, last paragraph]{es}, the authors 
express their belief that the Chow ring of invariant quotients of 
affine spaces, should be generated by Chern classes of `naturally 
given vector bundles'. Actually, they make this statement for the 
more restrictive case that they are considering, when the 
representation $G\rar\Gl(V)$ contains the homotheties. 

Following \cite[example 3, page 26]{sga6} we define the representation 
ring $\cal R(G)$ of $G$ to be the ring generated by isomorphism classes 
of representations of $G$ modulo the ideal generated by $[F]-[F']-[F'']$, 
where 
$$
0\lar F'\lar F\lar F''\lar 0
$$
is an exact sequence of $G$-modules; if ${\rm char\,}K=0$, reductive 
groups are linearly reductive, and in this case $F=F'\oplus F''$ as 
a $G$-module. The addition in $\cal R(G)$ is given by direct sum, and 
the product by tensor product. 

The restriction to the maximal torus defines a ring homomorphism 
$$
\cal R(G)\lar \cal R(T)^W,
$$
where on the right-hand-side we consider the $W$-invariant 
representations of $T$. By \cite[th\'eor\`eme 1.1, page 26]{sga6}, 
this homomorphism is actually an isomorphism. Since the representation 
ring of $T$ is isomorphic to the symmetric algebra of its group of 
characters, we deduce the isomorphism 
\begin{align}{\label{repres}}
\cal R(G)_{\mbb Q}
\srel{\cong}{\lar} 
\bigl(\Sym_{\mbb Q}^\bullet\cal X^*(T)_{\mbb Q}\bigr)^W
\cong 
(A^*_T)^W_{\mbb Q}\cong (A^*_G)_{\mbb Q}. 
\end{align}
Moreover, the composition $\cal R(G)_{\mbb Q}\rar (A^*_G)_{\mbb Q}$ 
is given by the Chern character 
$$
\cal R(G)\ni
(F,\rho)\lmt\ch 
\bigl(
\EG\times_\rho F\rar\BG
\bigr)
\in 
A^*_G.
$$
By abuse of notation, $\EG$ stands here for an appropriate open 
subset of an affine space, needed to define equivariant classes, 
and $\BG:=\EG/G$. 

\begin{proposition}{\label{prop:bdls}}
{\rm (i)} Let $\ell:={\rm rank\,}G=\dim T$.  There are at most 
$\ell$ elements in $\cal R(G)$ whose Chern classes generate 
the ring  $(A^*_T)^W_{\mbb Q}$. 

\nit{\rm (ii)} Consider $\chi\in\cal X^*(G)$ such that 
$\bbV^\sst(G,\chi)=\bbV^\s_{(0)}(G,\chi)$. The same statement 
as in {\rm (i)} holds for the Chow ring of the invariant quotient 
$\bbV\invq_\chi G$. 
\end{proposition}

\begin{proof}
Since $(A_*^T)^W_{\mbb Q}\rar A_*(Y)_{\mbb Q}$ is an epimorphism, 
it is enough to prove the first claim. Applying 
\cite[\S 5.5, th\'eor\`eme 4]{bou} to the $\mbb Q$-vector space 
$\cal X^*(T)_{\mbb Q}$, we deduce that there are algebraically 
independent, homogeneous elements 
$I_1,\dots,I_\ell\in (A_*^T)^W_{\mbb Q}$, such that 
$(A_*^T)^W_{\mbb Q}=\mbb Q[I_1,\dots,I_\ell]$. The isomorphism 
\eqref{repres} implies that there are elements 
$F_1,\dots,F_\ell\in\cal R(G)$ such that 
$$
\ch\,_{\deg I_j}(F_j)=I_j,\quad\forall\,j=1,\dots,\ell, 
$$
so that the Chern classes of $F_1,\dots,F_\ell$ generate 
$(A_*^T)^W_{\mbb Q}$. 
\end{proof}

A particularly comfortable situation arises when 
$G=\hbox{\Large$\times$}_{j=1}^{s}\Gl(n_j)$ is a product 
of linear groups; this happens, for instance, in the case 
of quiver representations. Then the Chow ring of the 
corresponding quotients will be generated by the Chern classes 
of the identical representations of the $s$ factors of $G$.


\section{The cohomology ring of the quotients}{\label{sct:cohom}}

In this section we assume that we work over the field of 
complex numbers, and moreover that the ring of invariants 
$\mbb C[\bbV]^G=\mbb C$. This latter condition guarantees 
that the invariant quotients $\bbV\invq_\chi G$, 
$\chi\in\cal X^*(G)$, are projective. 

\begin{theorem}{\label{thm:cohom}}
Let $\chi\in\cal X^*(G)$ be a character such that 
$\bbV^\sst(G,\chi)=\bbV^\s_{(0)}(G,\chi)$. Then the cohomology ring 
$$
H^*(\bbV\invq_\chi G;\mbb Q)\cong
{\bigl(H^*_T\bigr)}^W_{\mbb Q}
\biggl/
\wp_G\bigl\lan [\bbE(\l)]_T;\l\in\euf F(\chi)\bigr\ran_{\mbb Q}
\biggr.. 
$$
A fortiori, we deduce that the cycle map 
$$
{\rm cl}:A_*\bigl(\bbV\invq_\chi G\bigr)_{\mbb Q}
\lar 
H^*(\bbV\invq_\chi G;\mbb Q)
$$ 
is an isomorphism. 
\end{theorem}

\begin{proof}
The idea is to use the results of section \ref{sct:stab2} and to 
reduce our question to a similar one for actions of reductive 
groups on projective varieties, for which we know the cohomology 
of the semi-stable locus. 

We have proved in proposition \ref{prop:lin} that 
$$
\bbV^\sst(G,\chi)
=\mbb P(\mbb C\oplus V)^\sst\bigl(G,\Si_{(d,\chi^c)}\bigr)=:\O.
$$
For the convenience of the writing, we shall denote 
$\bar \bbV:=\mbb P(\mbb C\oplus V)$. 

Our first claim is that the restriction homomorphism 
$H^*_G(\bbV)\rar H^*_G(\O)$ is surjective: indeed, the diagram 
$$
\xymatrix{
H^*_G(\bar \bbV)
\ar[r]^-{\jmath^*_G}\ar[dr]_-{\alpha^*}
&
H^*_G(\bbV)
\ar[d]
\\ 
&
H^*_G(\O)
}
$$
commutes, and we know from \cite[page 175]{mfk} that $\alpha^*$ 
is surjective. Secondly, we claim that $\jmath^*_G$ is surjective 
too. This can be seen as follows: from the commutative diagram 
$$
\xymatrix{
\BG\cong\EG\times_G[1,0]\,
\ar@{^(->}[r]^{\kern4.5ex \jmath_0}\ar@{=}[dr]
&
\EG\times_G\bar \bbV
\ar[d]\ar@{}[drr]|{\kern-3em\Longrightarrow}
&
H^*_G
\ar[r]\ar@{=}[dr]_-{\rm id}
&
H^*_G(\bar \bbV)
\ar[d]^-{{(\jmath_0)}^*_G}
\\ 
&
\BG
&
&
H^*_G
}
$$
we deduce that ${(\jmath_0)}^*_G$ is surjective. As usual, 
$\EG\rar\BG$ stays for a universal $G$-bundle. The surjectivity 
of $\jmath^*_G$ is implied now by the commutative triangle 
\begin{align}{\label{eq:restr}}
\xymatrix{
H^*_G(\bar \bbV)
\ar[r]^-{\jmath^*_G}\ar@{->>}[dr]_-{{(\jmath_0)}^*_G}
&
H^*_G(\bbV)
\ar[d]^-\cong
\\
&
H^*_G
}
\end{align}
corresponding to the inclusions $\{0\}\subset \bbV\subset\bar \bbV$. 

The two claims together imply that 
\begin{align}{\label{eq:tg}}
\Ker
\bigl(
H^*_G(\bbV)\rar H^*_G(\O)
\bigr)
=
\jmath^*_G\Ker
\bigl(
H^*_G(\bar \bbV)\rar H^*_G(\O)
\bigr).
\end{align}
Now we are going to use the result \cite[corollaire 1.1]{br1} of 
M.~Brion, which expresses the right-hand-side in terms of the 
$W$-invariant part of the $T$-equivariant cohomology of the 
$T$-unstable locus. For this, we consider the diagram
$$
\xymatrix@C=1.5cm{
H^*_T(\bbV)
&
\kern1.5ex H^*_G(\bbV)
\ar@{>->}[l]_-{\vphi^*}
\\ 
H^*_T(\bar \bbV)
\ar[u]^-{\jmath^*_T}
&
\kern1.5ex H^*_G(\bar \bbV)
\ar@{>->}[l]^-{\vphi^*}\ar[u]_-{\jmath^*_G}
}
$$
Applying $\vphi^*$ to both sides of \eqref{eq:tg}, we obtain that 
$$
\begin{array}{ll}
\vphi^*\Ker
\bigl(
H^*_G(\bbV)\rar H^*_G(\O)
\bigr)
&
=
(\vphi^*\circ\jmath^*_G)\Ker
\bigl(
H^*_G(\bar \bbV)\rar H^*_G(\O)
\bigr)
\\[1.5ex]
& 
=
(\jmath^*_T\circ\vphi^*)
\Ker\bigl(
H^*_G(\bar \bbV)\rar H^*_G(\O)
\bigr)
\\[1.5ex]
&
\kern-1.8ex\srel{loc.\,cit.}{=}
\jmath^*_T\biggl[
\Ker\biggl(
H^*_T(\bar \bbV)\rar 
H^*_T\left(
\bar\bbV^\sst(T,\Si_{(d,\chi^c)})
\right)
\biggr)^{\!\lower.4ex\hbox{$_W$}}
\biggr]
\\[1.9ex]
&
=\biggl[
\jmath^*_T
\Ker\biggl(
H^*_T(\bar \bbV)\rar 
H^*_T\left(
\bar\bbV^\sst(T,\Si_{(d,\chi^c)})
\right)
\biggr)
\biggr]^{\lower.1ex\hbox{$_W$}}.
\end{array}
$$
The last equality holds because the projection onto the $W$-invariant 
part is a linear map (a Reynolds type homomorphism), which commutes 
with the pull-back. Let us point out that in order to apply Brion's 
result, we must be able to linearize the $G$-action in 
$\crl O_{\bar\bbV}(1)$; according to remark \ref{rmk:1}, this 
can be achieved after replacing $G$ with a suitable finite cover 
(see \ref{rmk:1}). This modification does not affect the cohomology 
ring, since we work with rational coefficients. 

The semi-stable locus $\bar\bbV^\sst(T,\Si_{(d,\chi^l)})$ corresponds 
to the representation 
$$
\bar\rho:T\lar \Gl(\mbb C\oplus V),\quad 
\bar\rho(g)=
\left(
\begin{array}{cc}
\chi(g)^{-1/d} & 0
\\[1ex] 
0 & \chi(g)^{-1/d}\rho(g)
\end{array}
\right),
$$
of the $d$-sheeted cover of $G$, for which the fractional power is 
well defined. 
If $\eta_1,\dots,\eta_r$ are the characters of the $T$-action on 
$\bbV\cong\mbb A^r_{\mbb C}$ (together with their multiplicities), 
then $T$ acts on $\mbb A^1_{\mbb C}\times\mbb A^r_{\mbb C}$ by the 
characters 
$$
\bar\eta_0:=\chi^{-1/d},\bar\eta_1:=\chi^{-1/d}\eta_1,\dots,
\bar\eta_r:=\chi^{-1/d}\eta_r.
$$ 
The irreducible components of the $T$-unstable locus in $\bar \bbV$ are 
of the form $\mbb P(F)$, where $F\hra\mbb C\oplus\mbb C^r$ runs over 
the set of unstable coordinate planes. For $F$ such a plane, we denote  
$z(F)\subset\{0,1,\dots,r\}$ its defining equations, that is 
$$
F=\{x_j=0\mid j\in z(F)\}\subset \mbb C\oplus\mbb C^r.
$$
Then \cite[th\'eor\`eme 2.1]{br1} says that 
$$
\Ker\biggl(
H^*_T(\bar \bbV)\rar 
H^*_T\left(
\bar\bbV^\sst(T,\Si_{(d,\chi^l)})
\right)
\biggr)
=
\biggl\lan
\prod_{j\in z(F)}(h+\bar\eta_j)
\ \biggl|\ 
\begin{array}{l}
F\subset\mbb C\oplus\mbb C^r\text{ is}
\\
\text{unstable plane}
\end{array}
\biggr\ran,
$$
where $h:=c_1\bigl(\crl O_{\bar \bbV}(1)\bigr)$. 

\unbar{\sl Claim}\quad  
$\jmath^*_T(h+\bar\eta_j)
=
\begin{cases}
\hskip0.3ex 0&\text{ if }0\in z(F);
\\ 
\eta_j&\text{ if }0\not\in z(F).
\end{cases}$

\begin{proof}(of claim) We notice first that 
$h+\bar\eta_j=c_1^T\bigl(\crl O(1)_{\bar\eta_j}\bigr)$, where we 
have denoted $\crl O(1)_{\bar\eta_j}$ the $T$-linearized invertible 
sheaf $\crl O_{\bar \bbV}(1)$ endowed with the $T$-action through 
the character $\bar\eta_j$. More precisely, the $T$-action on the 
geometric line bundle is  
$$
T\times\euf O_{\bar \bbV}(1)\lar\euf O_{\bar \bbV}(1),\quad 
t\times[(a,v),z]:=[t\times (a,v),\bar\eta_j(t)z]. 
$$
A consequence of the commutative diagram \eqref{eq:restr}, this 
time for $T$-equivariant cohomologies, is that for 
computing the restriction $\jmath^*_T(h+\bar\eta_j)\in H^2_T(\bbV)$, 
we must determine the character with which $T$ operates on 
$\jmath_0^*\euf O(1)_{\bar\eta_j}$. An immediate computation shows 
that it is 
$$
\bar\eta_0^{-1}\bar\eta_j=
\begin{cases}
\hskip 0.3ex1&\text{ if } j=0;
\\ 
\eta_j&\text{ if } j\neq 0.
\end{cases}
$$
This finishes the proof of the claim. 
\end{proof}

We continue now the proof of the theorem. First of all, we have proved 
in proposition \ref{prop:lin} that 
$\bbV^\sst(T,\chi)=\bbV\cap\bar\bbV^\sst\bigl(T,\Si_{(d,\chi^\ell)}\bigr)$. 
This implies that the intersections $\bbV\cap\mbb P(F)$, $F$ unstable 
plane in $\mbb C\oplus V$, coincide with the linear spaces 
$\bbE(\l)\hra\bbV$, $\l\in\euf F(\chi)$. For any such $F$ we distinguish 
between two possibilities: either $0\in z(F)$ or $0\not\in z(F)$. 
The claim above implies that 
$$
\jmath^*_T\biggl(
\prod_{j\in z(F)}(h+\bar\eta_j)
\biggr)
=
\begin{cases}
\hskip 2.3ex 0&\text{ if }0\in z(F); 
\\[1ex]\disp 
\prod_{j\in z(F)}\eta_j&\text{ if } j\neq 0.
\end{cases}
$$
Now, on one hand we know that $\bbV\cap \mbb P(F)$ coincides with some 
$\bbE(\l)$, $\l\in\euf F(\chi)$, and all the $\bbE(\l)$'s occur in this 
way, and on the other hand that the defining equations of 
$\bbV\cap\mbb P(F)$ are $\{x_j=0,j\in z(F)\}$. These two facts imply that 
$$
\prod_{j\in z(F)}\eta_j=[\bbE(\l)]_T,\text{ for appropriate }
\l\in\euf F(\chi),
$$
which finishes the proof of the theorem.
\end{proof}


\section{Construction of families of quotients}{\label{sct:fam}}

Let us recall that to any scheme $S$ and locally free sheaf 
$\crl F\rar S$, 
one can associate the Grassmannian ${\rm Grass}(d,\crl F)\rar S$ of 
$d$-dimensional quotients of $\crl E$. With this motivation in mind, 
we wish to address the problem of constructing families of varieties, 
which are invariant quotients of affine spaces, over arbitrary bases. 
Our construction relies on Seshadri's construction of quotients for 
actions of reductive group schemes. 

We fix a complex, connected, reductive group $G_{\mbb C}$ and a 
complex representation $\rho:G_{\mbb C}\rar\Gl(V_{\mbb C})$ of it 
with the property that $\mbb C[\bbV_{\mbb C}]^{G_{\mbb C}}=\mbb C$. 
Then $V_{\mbb C}$ decomposes uniquely, as a $G_{\mbb C}$-module, 
into its isotypical components 
$$
V_{\mbb C}=
\sum_{\omega\in\O}M_{\mbb C,\omega}\otimes_{\mbb C}V_{\mbb C,\omega},
\quad
\dim_{\mbb C}V_{\mbb C,\omega}=:\nu_\omega\geq 1,
$$ 
where the $M_{\mbb C,\omega}$ and the $V_{\mbb C,\omega}$'s are 
simple, and respectively trivial $G_{\mbb C}$-modules. We denote by 
$\rho_\omega:G_{\mbb C}\rar\Gl(M_{\mbb C,\omega})$ the corresponding 
representations. Further, for each $\omega$,  we choose 
an {\em admissible lattice} $M_\omega\subset M_{\mbb C,\omega}$ 
(see \cite[page 225]{chev}). Chevalley associates to this data, 
in \cite[section 4]{chev}, a reductive group scheme 
$\bG\rar\Spec\mbb Z$ which contains a maximal torus, and whose geometric 
fibres are connected, reductive groups having the same root system 
as $G_{\mbb C}$. Moreover, the representations $\rho_\omega$ extend to  
$$
\rho_\omega:\bG\lar\bGl(M_\omega):=
\Spec\biggl(\bigl(
\Sym_{\mbb Z}^\bullet\End_{\mbb Z}(M_\omega)\bigr)[\det^{-1}]\biggr). 
$$
We fix a character $\chi:\bG\rar\bG_m=\Spec\mbb Z[t,t^{-1}]$. 
Then for any field $K$, we get an induced character 
$\chi_K:\bG_K\rar\bG_{m,K}=\Spec K[t,t^{-1}]$. 

\begin{lemma}{\label{lm:triv-inv}}
If $G\rar\Gl(\nu;\mbb C)$ is a representation such that 
$\mbb C[t_1,\dots,t_\nu]^{G_{\mbb C}}=\mbb C$, then for 
any reduced ring $B$ holds $B[t_1,\dots,t_\nu]^{\bG_B}=B$. 
\end{lemma}

\begin{proof}
Assume the contrary, that there is $f\in B[t_1,\dots,t_\nu]^{\bG_B}$, 
homogeneous with $\deg f>0$. Since $B$ is reduced, there is 
$\mfrak q\in\Spec B$ such that the coefficients of $f$ are not contained 
in $\mfrak q$. We let $K$ to be the algebraic closure of quotient field 
$Q(B/\mfrak q)$, and we observe that the polynomial $F_K$ obtained by 
extending the coefficients of $f$ is non-zero, $\bG_K$-invariant, 
$\deg f_K>0$. Using \cite[theorem 1]{se}, we deduce that there 
is a homogeneous $F\in\mbb Z[t_1,\dots,t_\nu]^{\bG}$ with 
$\deg F>0$; this contradicts our hypothesis. 
\end{proof}

We consider further a ring $R$, which is a finite algebra over 
a universally Japanese ring; the examples we have in mind are 
$R=\mbb Z$ (for mixed characteristic methods) and  
$\mbb Z\bigl[[t]\bigr]$ (for the study of deformations). 

Let $S\rar\Spec R$ be a separated and reduced scheme of finite 
type. For  $\omega\in\O$, we consider a locally free sheaf 
$\crl V_\omega\rar S$ with $\rk_{\crl O_X}\crl V_\omega=\nu_\omega$, 
and define  
$$
\crl V:=\bigoplus_\omega M_\omega\otimes_{\mbb Z}\crl V_\omega\rar S. 
$$ 
We denote $\bV:=\uSpec\bigl(
\Sym^\bullet_{\crl O_S}\crl V\bigr)\srel{\pi}{\lar}S$ 
the corresponding geometric vector bundle. 

\begin{theorem}{\label{thm:family}} 
\nit{\rm (i)} There is a natural action 
$\bG\times_{\Spec\mbb Z}\bV\rar\bV$ of the  group scheme 
$\bG\rar\Spec\mbb Z$, which covers the trivial action on $S$. 

\nit{\rm (ii)} There is a $\bG$-invariant open subscheme 
$\bV^\sst(\bG,\chi)\subset\bV$ satisfying the following properties: 

${\rm (ii_a)}$ There is a categorical quotient $(\bY,q)$ for the 
$\bG$-action on $\bV^\sst(\bG,\chi)$. Moreover, there is a natural 
isomorphism 
$$
\bY\cong 
\uProj\biggl(\bigoplus_{n\geq 0}
\bigl(\pi_*\crl O_{\bV}\bigr)^{\bG}_{\chi^n}\biggr)
\lar S, 
$$
and $\bY$ is projective over $S$;
 
${\rm (ii_b)}$ For any algebraically closed field $K$, and any 
morphism $\Spec K\rar S$, 
$\Spec K\times_S\bV^\sst(\bG,\chi)$ consists of the 
$\chi_K$-semi-stable points of $\Spec K\times_S\bV$;
\end{theorem}

\begin{proof}
The proof will be divided in three steps: first we show that there 
is a natural $\bG$-action on $\bV$, next that the quotient exists 
locally, and finally we glue the local quotients together.  

For proving that there is a global $\bG$-action on $\bV$, is enough 
to prove that $\bGl(M_\omega)$ acts, for each $\omega\in\O$. We observe 
that the ring homomorphism 
$$
\begin{array}{rcl}
\Sym_{\crl O_S}^\bullet
\bigl(M_\omega\otimes_{\mbb Z}\crl V_\omega\bigr)
&
\lar
& 
\Sym_{\crl O_S}^\bullet
\bigl(\crl O_S\otimes_{\mbb Z}\End_{\mbb Z}(M_\omega)\bigr)
\otimes_{\crl O_S}
\Sym_{\crl O_S}^\bullet
\bigl(M_\omega\otimes_{\mbb Z}\crl V_\omega\bigr)
\\[1ex] 
m\otimes v
&
\lmt
&
\mu_\omega(m)\otimes v,
\end{array}
$$
with $\mu_\omega:\Sym_{\mbb Z}^\bullet M_\omega\rar
\Sym_{\mbb Z}^\bullet M_\omega\otimes \Sym_{\mbb Z}^\bullet M_\omega$ 
the co-multiplication, defines the $\bGl(M_\omega)$ action on 
$\uSpec(
\Sym_{\crl O_S}^\bullet M_\omega\otimes_{\mbb Z}\crl V_\omega)$. Since  
$$
\bV={\hbox{\huge$\times$}}_{\!S,\omega\in\O}
\uSpec(\Sym_{\crl O_S}^\bullet 
M_\omega\otimes_{\mbb Z}\crl V_\omega),
$$
$\bGl(M_\omega)$ acts on $\bV$ too. Repeating the very same 
construction as in section \ref{sct:stab}, we deduce also the 
existence of a $\bG$-action on $\bbA^1_S\times_S\bV$, where 
$\bG$ acts on $\bbA^1_S$ by $\chi$. 

What prevents us from constructing the quotient globally, is that 
$\bV$ can not be embedded as a $\bG$-invariant closed subset of an 
affine space over $\Spec R$; in this case, each $\crl V_\omega\rar S$ 
would be globally generated. What we do instead, is to make the 
constructions locally over $S$. 

Let us choose a finite covering ${(S_i=\Spec B_i)}_i$ of $S$ with 
open, affine subsets, which are moreover trivializing for all the 
$\crl V_\omega$'s. For every index $i$, we find then an integer 
$n_i\geq 1$ such that we have $S_i\hra \bbA_R^{n_i}$. It follows 
that 
$$
\begin{array}{ll}
\bV_i
& 
:={\hbox{\huge$\times$}}_{\!S,\omega\in\O}
\uSpec\bigl(
\Sym_{\crl O_S}^\bullet M_\omega\otimes_{\mbb Z}\crl V_\omega
\bigr)|_{S_i}
\\[1.5ex]
&\hra 
{\hbox{\huge$\times$}}_{\!\bbA_R^{n_i},\omega\in\O}
\uSpec\bigl(\Sym_{{\crl O}\hbox{$_{\bbA_R^{n_i}}$}}^\bullet 
M_\omega\otimes_{\mbb Z}\crl O_{\bbA_R^{n_i}}^{\oplus\nu_\omega}
\bigr)
=
\bbA_R^{n_i+\nu},
\end{array}
$$
with $\nu:=\sum_\omega\rk\,M_\omega\cdot \nu_\omega$, and moreover 
this embedding is $\bG$-equivariant. 

For $R_i:=\Gamma\bigl(\bbA^1_S\times_S\bV_i,
\crl O_{\bbA^1_S\times_S\bV_i}\bigr)\cong B_i[t,t_1,\dots,t_\nu]$, 
\cite[theorem 2]{se} implies that
$$
R_i^{\bG}=\bigoplus_{m\geq 0}t^m\cdot
B_i[t_1,\dots,t_\nu]^{\bG}_{\chi^m}\subset R_i
$$
is a finitely generated $B_i$-algebra. We define the $\chi$-unstable 
locus $\bV_i^\us(\bG,\chi)\hra\bV_i$ to be the closed subscheme 
defined by the ideal
$$
\bigl\lan B_i[t_1,\dots,t_\nu]^{\bG}_{\chi^m}\,;\,m\geq 1\bigr\ran
\subset
B_i[t_1,\dots,t_\nu]
$$
and let $\bV_i^\sst(\bG,\chi)$ to be its complement. The glueing 
argument of \cite[theorem 1.10]{mfk} together with 
\cite[theorem 3 and remark 8]{se} imply that the categorical 
quotient $\bY_i$ of $\bV_i^\sst(\bG,\chi)$ by $\bG$ exists, 
and it is isomorphic to $\uProj\bigl(
\bigoplus_{m\geq 0}B_i[t_1,\dots,t_\nu]^{\bG}_{\chi^m}
\bigr)$. It follows now from \cite[proposition 5.5.1]{ega2} that 
this quotient is projective over 
$\uSpec\bigl(B_i[t_1,\dots,t_\nu]^{\bG}\bigr)$. 

Using \cite[lemma 2]{se}, we remark that for any two open subsets 
$S_i,S_j\subset S$, holds 
$$
\bV_i^\sst(\bG,\chi)|_{S_i\cap S_j}
=
\bV_j^\sst(\bG,\chi)|_{S_i\cap S_j},
$$
so that there is a well-defined open subset $\bV^\sst(\bG,\chi)
\subset\bV$ of $\chi$-semi-stable points. Finally, using the 
universality property of categorical quotients, we glue 
together the $\bV_i^\sst(\bG,\chi)\rar\bY_i$'s into a scheme 
$\bY$, and this comes with a natural morphism 
$\bV^\sst(\bG,\chi)\rar\bY$. 

It remains to prove the behaviour of the semi-stable locus under 
base change; the question being local on $S$, we may assume that 
$S=\Spec B$. The inclusion 
$K\times_B \bV^\sst(G,\chi)\subset\bV_K^\sst(\bG_K,\chi_K)$ is 
trivial, so that we prove the converse. Consider a $K$-valued 
point $x_0\in\bV^\sst(\bG_K,\chi_K)$, and a homogeneous 
$f\in K[t_1,\dots,t_\nu]^{\bG_K}_{\chi_K^m}$, $m\geq 1$, such that 
$f(x_0)\neq 0$. 
Then $F_K:=t^{\deg f}f\in K[t,t_1,\dots,t_\nu]^{\bG_K}$, and 
applying \cite[theorem 1]{se} we find a homogeneous 
$F\in B[t,t_1,\dots,t_\nu]^{\bG_B}$ such that $\deg F>0$ and 
$F(1,x_0)\neq 0$. By lemma \ref{lm:triv-inv}, 
$B[t_1,\dots,t_\nu]^{\bG_B}=B$, so that 
$F\in\bigoplus_{m\geq 1}t^m B[t_1,\dots,t_\nu]^{\bG_B}_{\chi^m}$, 
and therefore $x_0\in K\times_B \bV^\sst(G,\chi)$. 
\end{proof}

\begin{remark}{\label{rmk:hope}}
In case $S=\Spec\mbb Z$, we get a {\em flat} family of quotients 
parameterized by $\Spec\mbb Z$, and one might hope to relate the 
complex geometric properties of the quotients to the arithmetical 
ones. 
\end{remark}

As we have already mentioned, one of the most common situation 
where the construction applies is that of the Grassmann bundle 
of $d$-dimensional quotients of a locally free sheaf 
$\crl F\rar S$. In this case, we define 
$\crl V:=\Hom_{\mbb Z}(\mbb Z^d,\crl F)\cong 
\mbb Z^d\otimes_{\mbb Z}\crl F$, $\bG:=\bGl(d)$, $\chi:=\det$. 
\medskip

As a first application, we describe how does the cone decomposition 
obtained in section \ref{sct:fan} varies with the characteristic of 
the ground field. 

\begin{corollary}{\label{cor:cone}}
Assume $S=\Spec\mbb Z$, 
$\bV=\Spec\bigl(\Sym^\bullet_{\mbb Z}M\bigr)$, and 
$\chi,\veps:\bG\rar\bG_m$ are characters such that 
$\veps_{\mbb C}\in C(\chi_{\mbb C})$. Then there is a prime 
$p_0>0$ depending on $\veps$, such that $\veps_K\in C(\chi_K)$ 
for any algebraically closed field $K$ of characteristic zero 
or ${\rm char\,}K>p_0$. 
\end{corollary}

\begin{proof}
The first remark is that 
$\bV_{\mbb Q}^\sst(\bG_{\mbb Q},\chi_{\mbb Q})
\subset
\bV_{\mbb Q}^\sst(\bG_{\mbb Q},\veps_{\mbb Q})$: indeed, otherwise 
there would exist 
$\mfrak q\in\Spec\bigl(\Sym^\bullet_{\mbb Q}M_{\mbb Q}\bigr)$ 
such that 
$$
\mfrak q\in\bV_{\mbb Q}^\sst(\bG_{\mbb Q},\chi_{\mbb Q}) 
\quad\text{and}\quad 
\mfrak q\supset\bigl\lan (\Sym^\bullet_{\mbb Q}M_{\mbb Q})
^{\bG_{\mbb Q}}_{\veps_{\mbb Q}^m}\,;\,m\geq 1\bigr\ran.
$$ 
Using theorem \ref{thm:family} $\rm (ii)_b$, we obtain that 
$\mfrak q_{\mbb C}\in\bV_{\mbb C}^\sst(\bG_{\mbb C},\chi_{\mbb C})$ 
and 
$$
\mfrak q_{\mbb C}\supset
\sqrt{
\bigl\lan (\Sym^\bullet_{\mbb Q}M_{\mbb Q})
^{\bG_{\mbb Q}}_{\veps_{\mbb Q}^m}\,;\,m\geq 1\bigr\ran_{\mbb C}
}
=
\sqrt{
\bigl\lan (\Sym^\bullet_{\mbb C}M_{\mbb C})
^{\bG_{\mbb C}}_{\veps_{\mbb C}^m}\,;\,m\geq 1\bigr\ran
}\,,
$$
which contradicts our hypothesis. 

Applying theorem \ref{thm:family} $\rm (ii)_b$ again, we deduce 
the validity of the corollary in characteristic zero. To settle 
the problem in positive characteristic, let us notice that the 
irreducible components of $\bV^\us(\bG,\veps)_{\rm red}\hra\bV$ 
are divided in two groups, according to whether their generic 
point is mapped onto $\lan 0\ran\in\Spec\mbb Z$ or not ($i.e.$ 
the residue field has characteristic zero or strictly positive). 
Let $U\subset\Spec\mbb Z$ be the complement of the finite set 
of primes where we have such `bad reduction'. Since 
$\bV^\us(\bG_{\mbb Q},\veps_{\mbb Q})\hra 
\bV^\us(\bG_{\mbb Q},\chi_{\mbb Q})$, we deduce that 
$\bV^\us(\bG,\veps)|_U\hra \bV^\us(\bG,\chi)|_U$. 
\end{proof}

We specialize now to the case where $S$ is a reduced and irreducible 
$K$-scheme of finite type, where $K$ is an algebraically closed 
field. Our construction yields in this case a locally trivial fibration 
$\bY\rar S$, whose fibres are isomorphic to $\mbb V\invq_\chi G$, 
and we interested in computing its Chow ring. 

Reasoning locally over $S$, we see that the 
$\chi|_{\lower.3ex\hbox{$_T$}}$-unstable locus of $\bV$ for the 
induced $T$-action is the union 
$$
\bV^\us(T,\chi)=\bigcup_{\l\in\euf F(\chi)}\bE(\l), 
$$
where $\bE(\l):=\uSpec\bigl(\crl V/\crl E(\l)\bigr)\hra\bV$, 
$\l\in\euf F(\chi)$, are sub-vector bundles associated to 
locally free subsheaves $\crl E(\l)\subset\crl V$. 
The $\chi$-unstable locus for the $G$-action is then 
$$
\bV^\us(G,\chi)=\bigcup_{\l\in\euf F(\chi)}G\cdot\bE(\l). 
$$
We still denote by 
$\wp_G:A_*(S)_{\mbb Q}\otimes_{\mbb Q} (A_*^T)_{\mbb Q}
\rar A_*(S)_{\mbb Q}\otimes_{\mbb Q} (A_*^T)_{\mbb Q}^W$ the 
homomorphism obtained from $\wp_G$, as defined in section 
\ref{sct:chow}, by extending the scalars. 

\begin{theorem}{\label{thm:chow-family}}
Assume that the character $\chi\in\cal X^*(G)$ is such that 
$\bbV^\sst(G,\chi)=\bbV^\s_{(0)}(G,\chi)$. Then: 

\nit{\rm (i)} the quotient $\bY:=\bV^\sst(G,\chi)/G$ is geometric; 

\nit{\rm (ii)} the Chow ring 
$$
A_*(\bY)_{\mbb Q}
\cong
A_*(S)_{\mbb Q}\otimes_{\mbb Q} (A_*^T)_{\mbb Q}^W
\biggl/
\wp_G
\bigl\lan
\sfe^T(\crl E(\l))\,;\,\l\in\euf F(\chi)
\bigr\ran_{\mbb Q}.
\biggr.
$$

\nit{\rm (iii)} Assume moreover that $K=\mbb C$. Then 
$$
H^*(\bY;\mbb Q)
\cong
H^*(S)_{\mbb Q}\otimes_{\mbb Q} (H^*_T)_{\mbb Q}^W
\biggl/
\wp_G
\bigl\lan
\sfe^T(\crl E(\l))\,;\,\l\in\euf F(\chi)
\bigr\ran_{\mbb Q}.
\biggr.
$$ 
\end{theorem}

\begin{proof}
The quotient $\bV^\sst(G,\chi)\rar\bY$ is geometric because it is so 
locally over $S$. Now we prove the second statement: since $G$ acts 
with finite stabilizers, we have the exact sequence
$$
A_*^G
\biggl(
\bigcup_{\l\in\euf F(\chi)}G\cdot\bE(\l)
\biggr)_{\mbb Q}
\lar
A_*^G(\bV)_{\mbb Q}
\lar 
A_*^G
\bigl(\bV^\sst(G,\chi)\bigr)_{\mbb Q}
\cong
A_{*-\dim G}(\bY)_{\mbb Q}
\lar 0. 
$$
Repeating the proof of proposition \ref{prop:image} in our relative 
setting, and using the excess intersection formula 
\cite[example 6.3.5]{fu}, we find that the image of the homomorphism 
$A_*^G\bigl(\bE(\l)\bigr)_{\mbb Q}\rar A_*^G(\bV)_{\mbb Q}$ is the 
ideal $\wp_G\bigl\lan\sfe^T\bigl(\crl E(\l)\bigr)\bigr\ran_{\mbb Q}$. 
This proves the second claim. 

The proof of the third statement is divided in three steps: 

-- We recall that $\bY\rar S$ is a locally trivial fibre bundle, 
with fibre $Y:=\bbV\invq_\chi G$. According to proposition 
\ref{prop:bdls}, the Chow ring of $Y$, hence its cohomology ring 
too, is generated by characteristic classes of (virtual) 
representations of $G$. In the relative setting, one may still 
consider the locally free sheaves associated to representations 
of $G$, and the corresponding characteristic classes yield a 
cohomology extension of the fibre. We conclude that the Leray-Hirsch 
theorem applies to $\bY\rar S$, and that the long exact sequence 
for the cohomology of the pair 
$\bigl(\bV,\bV^\sst(\bbV,\chi)\bigr)$ splits into short exact 
sequences 
$$
0\lar 
H^*_G\bigl(\bV,\bV^\sst(G,\chi);\mbb Q\bigr)
\srel{\imath_S^*}{\lar}
H^*_G(\bV;\mbb Q)
\srel{\jmath_S^*}{\lar} 
H^*_G\bigl(\bV^\sst(G,\chi);\mbb Q\bigr)
\lar 0.
$$

-- Since the ring homomorphism 
$H^*_G(\bbV;\mbb Q)\rar H^*_G\bigl(\bbV^\sst(G,\chi);\mbb Q\bigr)
\cong H^*(Y)$ is surjective, the short sequence
$$
0\lar 
H^*_G\bigl(\bbV,\bbV^\sst(G,\chi);\mbb Q\bigr)
\srel{\imath^*}{\lar}
H^*_G(\bbV;\mbb Q)
\srel{\,\jmath^*}{\lar} 
H^*_G\bigl(\bbV^\sst(G,\chi);\mbb Q\bigr)
\lar 0,
$$
is exact, and moreover that 
$\imath^*H^*_G\bigl(\bbV,\bbV^\sst(G,\chi);\mbb Q\bigr)
=\wp_G\bigl\lan [\bbE(\l)]_T\,;\,\l\in\euf F(\chi)\bigr\ran_{\mbb Q}$. 

-- The restriction to the fibres of the fibre bundle pair 
$\bigl(\bV,\bV^\sst(G,\chi)\bigr)\srel{\pi}{\lar}S$ of the ideal 
$\wp_G\bigl\lan\sfe^T\bigl(\crl E(\l)\bigr)\bigr\ran_{\mbb Q}
\subset H^*_G(\bV;\mbb Q)$, $\l\in\euf F(\chi)$, equals 
$\wp_G\bigl\lan [\bbE(\l)]_T\bigr\ran$. Moreover, we remark that 
the restriction of $\pi^*\crl E(\l)^\vee$ to $\bV^\sst(G,\chi)$ 
has a nowhere vanishing, $T$-equivariant section: it is the 
composition of the canonical section $\bV\rar\pi^*\crl V^\vee$ 
with the projection $\crl V^\vee\rar\crl E(\l)^\vee$. It follows 
that 
$$
\jmath_S^*\sfe^T\bigl(\crl E(\l)^\vee\bigr)
=(-1)^{{\rm rank}\crl E(\l)}
\jmath_S^*\sfe^T\bigl(\crl E(\l)^\vee\bigr)=0,
$$
and therefore $\jmath_S^*\bigl[
\wp_G\bigl\lan\sfe^T\bigl(\crl E(\l)\bigr)\bigr\ran_{\mbb Q}
\bigr]=0$. 
This implies in turn that the Leray-Hirsch theorem applies to 
the pair $\bigl(\bV,\bV^\sst(\bbV,\chi)\bigr)\rar S$, and also 
that 
$$
\imath_S^*\bigl[
H^*_G\bigl(\bV,\bV^\sst(G,\chi);\mbb Q\bigr)
\bigr]
=
\wp_G\bigl\lan
\sfe^T\bigl(\crl E(\l)\bigr)\,;\,\l\in\euf F(\chi)
\bigr\ran_{\mbb Q}.
$$
Since $\jmath_S^*$ is a ring homomorphisms, our claim follows. 
\end{proof}


\section{The Picard group, ample cone and the cohomology of line bundles}
{\label{coh-line-bdls}}

We wish to use the results obtained so far to compute the Picard group 
and the ample cone of the geometric quotients that we obtain. Secondly, 
we will be concerned about the vanishing of the higher cohomology groups 
of ample, and more generally nef, invertible sheaves. 

\begin{definition}{\label{line-bdls}}
Let $\chi\in\cal X^*(G)$ be a character, and consider the categorical 
quotient $q:\bbV^\sst(G,\chi)\rar Y_\chi:=\bbV^\sst(G,\chi)/G$. For 
a further character $\veps\in\cal X^*(G)$, we define 
$$
\crl L_{Y_\chi,\veps}:=
\bigl(q_*\crl O_{\bbV^\sst(G,\chi)}\bigr)^G_{\veps}.
$$
In other words, for each open subset $U\subset Y_\chi$, 
$\Gamma(U,\crl L_\veps)=
\Gamma\bigl(q^{-1}U,\crl O_{\bbV^\sst(G,\chi)}\bigr)^G_\veps$. 
\end{definition}

\nit It is easy to see that $\crl L_{Y_\chi,\veps}$ is torsion-free 
and coherent. If $\bbV^\sst(G,\chi)=\bbV^\s_{(0)}(G,\chi)$,  
$\crl L_{Y_\chi,\veps}$ has rank one, and if moreover the $G$-action  
on $\bbV^\sst(G,\chi)$ is free, then $\crl L_{Y_\chi,\veps}$ is 
actually invertible. 

We recall from section \ref{sct:fan} that for $\chi\in\cal X^*(G)$, 
$C(\chi)\in\Delta^G(\bbV)$ denotes the cone which contains $\chi$ 
in its relative interior. 

For a character $\chi\in\cal X^*(G)$ such that 
$\bbV^\sst(G,\chi)=\bbV^\s_{(0)}(G,\chi)\neq\emptyset$, we define 
$$
\begin{array}{rl}
\euf F'(\chi):=
&
\{\l\in\euf F(\chi)\mid\codim_\bbV\ G\cdot\bbE(\l)=1\}
\\[1ex]
=
&
\{\l\in\euf F(\chi)\mid\codim_\bbV\bbE(\l)+{\rm Stab}_G\bbE(\l)
=\dim G+1\}.
\end{array}
$$
Denoting $L_\l$ the Levi component of ${\rm Stab}_G\bbE(\l)$, 
we observe that for $\l\in\euf F'(\chi)$, 
$\deg[\bbE(\l)]_T-\deg\Delta_G/\Delta_{L_\l}=1$ (see section 
\ref{sct:chow} for the notations). Since in this case 
$\wp_G\lan[\bbE(\l)]_T\ran$ is a principal ideal, formula 
\eqref{local} implies that 
$$
\wp_G\lan[\bbE(\l)]_T\ran_{\mbb Q}
=\lan \wp_G[\bbE(\l)]_T\ran_{\mbb Q}.
$$ 
For shorthand, we will denote $\veps(\l):=\wp_G[\bbE(\l)]_T\in 
(A^1_T)^W=\cal X^*(T)^W\cong\cal X^*(G)$. 

\begin{proposition}{\label{cor:pic}}
Let $\chi\in\cal X^*(G)$ be a character such that 
$\bbV^\sst(G,\chi)=\bbV^\s_{(0)}(G,\chi)$. Then:  

\nit{\rm (i)} $\Pic(Y_\chi)_{\mbb Q}
\cong \cal X^*(G)_{\mbb Q}\bigl/\bigl\lan\veps(\l)\,;\,
\l\in\euf F'(\chi)\bigr\ran_{\mbb Q}$;

\nit{\rm (ii)} The ample cone of $Y_\chi$ is 
$$
\Pic(Y_\chi)^{\rm ample}_{\mbb Q}
\cong 
{\rm int.}\,\hat C(\chi)_{\mbb Q},
$$ 
where $\hat C(\chi)_{\mbb Q}:={\rm Image}\bigl[
\cal X^*(G)_{\mbb Q}\cap C(\chi)
\rar 
\cal X^*(G)_{\mbb Q}\bigr/\bigl\lan\veps(\l)\,;\,
\l\in\euf F'(\chi)\bigr\ran_{\mbb Q}
\bigr]$.
\end{proposition}

\begin{proof}
(i) Indeed,
$$
\begin{array}{ll}
{\Pic(Y_\chi)}_{\mbb Q}
&
=
{A^1(Y_\chi)}_{\mbb Q}
=
{(A^1_T)}^W_{\mbb Q}\biggl/ 
{(A^1_T)}^W_{\mbb Q}\cap 
\wp_G\bigl\lan [\bbE(\l)]_T\,;\,\l\in\euf F(\chi)\bigr\ran_{\mbb Q}
\biggr.
\\ 
&
=
\cal X^*(G)_{\mbb Q}\bigl/ 
\bigl\lan \veps(\l)\,;\,\l\in\euf F'(\chi)\bigr\ran_{\mbb Q}
\biggr..
\end{array}
$$

\nit(ii) We recall that any character 
$\veps\in\cal X^*(G)_{\mbb Q}\cap{\rm int.}C(\chi)$ defines the same 
invariant quotient as $\chi$, that is $Y$, but endows it with a possibly 
different polarization. Namely, for $n\geq 1$ sufficiently large, the  
sheaf $\crl L_{\veps^n}\rar Y$ is invertible and ample. 

For $\l\in\euf F'(\chi)$, let $f_\l\in K[\bbV]$ be the equation of 
$G\cdot\bbE(\l)\hra\bbV$; $G$ acts on $\mbb Qf_\l$ by a character, 
which is just the class $[G\cdot\bbE(\l)]_G\in A^*_G$, that is 
$\veps(\l)$. We conclude that $f_\l\in K[\bbV]^G_{\veps(\l)}$. 
We deduce further that for 
$\veps'=\sum_{\l\in\euf F'(\chi)} k_\l\veps(\l)$, with $k_\l\geq 0$, 
$$
\prod_{\l\in\euf F'(\chi)}\kern-1ex f_\l^{k_\l}:
\crl L_\veps
=\bigl(q_*\crl O_{\bbV^\sst(G,\chi)}\bigr)^{G,\veps}
\lar 
\bigl(q_*\crl O_{\bbV^\sst(G,\chi)}\bigr)^{G,\veps\veps'}
=\crl L_{\veps\veps'}
$$
is an isomorphism. The two remarks together prove the inclusion of 
the right hand side into the left hand side. 

For the converse inclusion, let us consider an ample line bundle 
$\crl L\rar Y$. It follows from the first part that 
$\crl L\cong\crl L_\veps$ for some $\veps\in\cal X^*(G)$. 
According to \cite[th\'eor\`eme 4.5.2]{ega2}, there is some 
$n\geq 1$ such that the non-vanishing loci of the sections in 
$\crl L^{\otimes n}$ cover $Y$ with affine open sets. 
This means that for any $x\in\bbV^\sst(G,\chi)$, there is 
$f\in K[\bbV^\sst(G,\chi)]^G_{\veps^n}$ such that $f(x)\neq 0$. 
We choose finitely many such $f$'s such that the non-vanishing 
loci cover $\bbV^\sst(G,\chi)$. The components of codimension one 
of $\bbV\sm\bbV^\sst(G,\chi)$ are the $G\cdot\bbE(\l)$'s, with 
$\l\in\euf F'(\chi)$. Therefore we find a monomial 
$\phi:=\prod_{\l\in\euf F'(\chi)}\kern-0.5ex f_\l^{k_\l}$ with 
the property that the $\phi f$'s, with $f$ as above, extend to 
regular functions on $\bbV$; they all belong to 
$K[\bbV]^G_{\veps^n\cdot\prod_\l\veps_\l^{k_\l}}$, and their 
non-vanishing locus contain $\bbV^\sst(G,\chi)$. This means that 
$\bbV^\sst(G,\chi)\subset
\bbV^\sst(G,\veps^n\cdot\prod_\l\veps_\l^{k_\l})$, 
and therefore $\veps^n\cdot\prod_\l\veps_\l^{k_\l}\in C(\chi)$. 

All together, proves that 
$$
{\rm int.}\,\hat C(\chi)_{\mbb Q}
\subset 
\Pic(Y_\chi)^{\rm ample}_{\mbb Q}
\subset
\hat C(\chi)_{\mbb Q}
\subset
\Pic(Y_\chi)_{\mbb Q}.
$$
The conclusion that 
$\Pic(Y_\chi)^{\rm ample}_{\mbb Q}
={\rm int.}\,\hat C(\chi)_{\mbb Q}$ follows, because we know from 
\cite[corollaire 4.5.8]{ega2} that the ample cone is open in the 
Picard variety.  
\end{proof}

\begin{remark}{\label{rmk:inj}}
If and $\codim_\bbV\bbV^\us(G,\chi)\geq 2$, there is a short exact 
sequence 
$$
0\lar\Pic(Y_\chi)\srel{q^*}{\lar} \cal X^*(G)\lar \cal T\lar 0,
$$
where $\cal T$ is by definition the cokernel of $q^*$. Indeed, the 
pull-back by $q$ of an invertible sheaf on $Y_\chi$ gives rise to 
a $G$-linearized invertible sheaf on $\bbV^\sst(G,\chi)$. The 
codimension condition ensures that this $G$-linearized invertible 
sheaf uniquely extends to a $G$-linearized sheaf on $\bbV$, which 
corresponds to a character of $G$. 

If moreover $\bbV^\sst(G,\chi)=\bbV^\s_{(0)}(G,\chi)$, the previous 
proposition says that $\cal T$ is a $\mbb Z$-torsion module, 
because a suitably large multiple of any character defines an 
invertible sheaf on $Y_\chi$. 
\end{remark}

Now that we have described the line bundles on our quotient varieties, 
we wish to compute their cohomology groups. It is well-known that, 
over complete toric varieties, the higher cohomology groups of nef 
line bundles vanish. What we are going to prove is that a similar 
result holds in our `non-abelian setting'. We must say from the very 
beginning that our result is a more or less straightforward consequence 
of the Hochster-Roberts theorem. 

First of all, we remark that for a pair of characters 
$\chi,\veps\in\cal X^*(G)$, with $\veps\in C(\chi)$, the natural 
inclusion $\bbV^\sst(G,\chi)\subset\bbV^\sst(G,\veps)$ induces 
a morphism on the level of quotients
\begin{align}{\label{cd}}
\xymatrix{
\bbV^\sst(G,\chi)
\ar[r]\ar[d]_-{q_\chi}
&
\bbV^\sst(G,\veps)
\ar[d]^-{q_\veps}
\\ 
Y_\chi
\ar[r]^-{\psi}
&
Y_\veps
}
\end{align}
which is projective and open. In particular, $\psi$ is surjective. 

\begin{lemma}{\label{lm:pre}} 
{\rm (i)} For $f\in K[\mbb V]^G_{\veps^m}$, with $m\geq 1$, the 
preimage under $\psi$ of the affine, open subset 
$D(f,Y_\veps)\subset Y_\veps$ defined by $f$ is isomorphic to 
$\Proj\bigl(\underset{n\geq 0}{\oplus}(f^{-1}K[\bbV])^G_{\chi^n}\bigr)$. 

\nit{\rm (ii)} Assume ${\rm char\,}K=0$. Then 
$$
\rR^i\psi_*\crl O_{Y_\chi}=0,\quad\forall\,i>0.
$$
\end{lemma}

\begin{proof} 
(i) If $D(f)\subset\bbV$ denotes the non-vanishing locus 
of $f$ in $\bbV$, then clearly 
$$
\psi^{-1}D(f,Y_\veps)=
\bigl(D(f)\cap\bbV^\sst(G,\chi)\bigr)\bigl/G\bigr.
=D(f)^\sst(G,\chi)\bigl/G\bigr..
$$ 
As $D(f)=\Spec\bigl(f^{-1}K[\mbb V]\bigr)$, we conclude using 
the same trick as in lemma \ref{lm:us}. 

\nit (ii) Is enough to prove that 
$\bigl(\rR^i\psi_*\crl O_{Y_\chi}\bigr)\bigl(D(f,Y_\veps)\bigr)
=H^i\bigl(\psi^{-1}D(f,Y_\veps),\crl O_{Y_\chi}\bigr)=0$, 
for any $f\in K[\mbb V]^G_{\veps^m}$, $m\geq 1$. We have just 
seen that $\psi^{-1}D(f,Y_\veps)=\Proj R$, with 
$$
R=\oplus_{n\geq 0}R_n
:=\bigoplus_{n\geq 0}(f^{-1}K[\bbV])^G_{\chi^n}
\text{ and also } 
\bigoplus_{n\geq 0}
(f^{-1}K[\bbV])^G_{\chi^n}=K[\bbA^1_K\times D(f)]^G. 
$$  
Since $G$ is linearly reductive, as ${\rm char\,}K=0$,and 
$\bbA^1_K\times D(f)$ is a regular, affine variety, we conclude 
that we are in the situation of \cite[theorem 3.4]{ke1}, and 
therefore the positively graded part of the local cohomology 
groups $\bigl[H^i_{R_+}(R)\bigr]_{\geq 0}$ vanish, for all 
$i\geq 0$, where $R_+:=\oplus_{n\geq 1}R_n$. The 
Grothendieck-Serre correspondence \cite[theorem 20.4.4]{bs} 
says that 
$$
\biggl[H^{i+1}_{R_+}(R)\biggr]_{\geq 0}
\cong
\bigoplus_{j\geq 0} H^i\bigl(\Proj R,\crl O_{\Proj R}(j)\bigr),
\quad\forall\,i\geq 1,
$$
and the vanishing of the degree zero part implies  
$H^i\bigl(\psi^{-1}D(f,Y_\veps),\crl O_{Y_\chi}\bigr)=0$.  
\end{proof}

\begin{lemma}{\label{glob-gen}}
Let $\chi\in\cal X^*(G)$ such that $\codim_\bbV\bbV^\us(G,\chi)\geq 2$, 
and $\veps\in\cal X^*(G)\cap C(\chi)$. Then the following hold: 

\nit{\rm (i)} $\psi_*\crl L_{Y_\chi,\veps}=\crl L_{Y_\veps,\veps}$;

\nit{\rm (ii)} if $\crl L_{Y_\veps,\veps}$ is invertible, then the 
natural homomorphism 
$$
\psi^*\crl L_{Y_\veps,\veps}
=\psi^*\psi_*\crl L_{Y_\chi,\veps}\lar\crl L_{Y_\chi,\veps}
$$ 
is an isomorphism. 
\end{lemma}

We remark that for any $\veps\in\cal X^*(G)$, some positive multiple 
of it fulfills the requirement of (ii) above. 

\begin{proof}
(i) For an open subset $U'\subset Y_\veps$ holds 
$$
\begin{array}{ll}
\Gamma\bigl(U',\psi_*\crl L_{Y_\chi,\veps}\bigr)
&
=
\Gamma\bigl(\psi^{-1}U',\crl L_{Y_\chi,\veps}\bigr)
=
\Gamma\bigl((\psi\circ q_\chi)^{-1}U',\crl O_\bbV\bigr)^G_\veps
\\[1ex] 
&
=
\Gamma\bigl(
q_\veps^{-1}U'\!\cap\!\bbV^\sst(G,\chi),\crl O_\bbV
\bigr)^G_\veps
\srel{(*)}{=}
\Gamma\bigl(q_\veps^{-1}U'\!,\crl O_\bbV\bigr)^G_\veps
=\Gamma\bigl(U'\!,\crl L_{Y_\veps,\veps}\bigr), 
\end{array}
$$
where $(*)$ holds because of the assumption that 
$\codim_\bbV\bbV^\us(G,\chi)\geq 2$. 

\nit (ii) Since $\crl L_{Y_\veps,\veps}$ is locally free, its 
pull-back is the same. Using remark \ref{rmk:inj}, we deduce 
that $\psi^*\crl L_{Y_\veps,\veps}\cong \crl L_{Y_\veps,\veps'}$, 
for some character $\veps'$. Consider now an affine, open subset 
$U'\subset Y_\veps$ over which $\crl L_{Y_\veps,\veps}$ is 
trivialized by a section 
$s'\in\Gamma(q_\veps^{-1}U',\crl O_\bbV)^G_\veps$. Then, on one 
hand, we have 
$$
\begin{array}{ll}
\Gamma\bigl(\psi^{-1}U',\psi^*\crl L_{Y_\veps,\veps}\bigr)
&
=
\Gamma\bigl(\psi^{-1}U',\crl L_{Y_\chi,\veps'}\bigr)
=\Gamma\bigl(
\bbV^\sst(G,\chi)\cap q_\veps^{-1}U',\crl O_\bbV
\bigr)^G_{\veps'}
\\[1ex] 
&
=
\Gamma\bigl(q_\veps^{-1}U',\crl O_\bbV\bigr)^G_{\veps'}.
\end{array}
$$
For the last step we use that $\codim_\bbV\bbV^\us(G,\chi)\geq 2$. 
On the other hand, 
$$
\begin{array}{ll}
\Gamma\bigl(\psi^{-1}U',\psi^*\crl L_{Y_\veps,\veps}\bigr)
&=
s'\cdot\Gamma\bigl(\psi^{-1}U',\crl O_{Y_\chi}\bigr)
=
s'\cdot\Gamma\bigl(
\bbV^\sst(G,\chi)\cap q_\veps^{-1}U',\crl O_{\bbV}
\bigr)^G
\\[1ex]
&
=s'\cdot\Gamma\bigl(q_\veps^{-1}U',\crl O_{\bbV}\bigr)^G
=\Gamma\bigl(q_\veps^{-1}U',\crl O_{\bbV}\bigr)^G_\veps,
\end{array}
$$
since $s'$ trivializes $\crl L_{Y_\veps,\veps}$ over $U'$. 
The two equalities imply $\veps'=\veps$. 
\end{proof}

We fix $\bG\rar\Spec\mbb Z$ a reductive group scheme, which contains 
a maximal torus, and whose geometric fibres are connected, reductive 
groups, all having the same root system; we consider a 
representation $\bG\rar\bGl(M)$, as constructed at the beginning 
of the previous section, and assume that 
$\bigl(\Sym_{\mbb Z}^\bullet M\bigr)^\bG=\mbb Z$. We define   
$\bV:=\Spec\bigl(\Sym_{\mbb Z}^\bullet M\bigr)$, and, for a field 
$K$, we let $\bbV_K:=\bV\times_{\Spec\mbb Z}\Spec K
=\Spec\bigl(\Sym_{K}^\bullet M\otimes_{\mbb Z}K\bigr)$. Further, 
we fix a character $\chi:\bG\rar\bG_m$ such that 
$\codim_{\bbV_{\mbb C}}
\mbb V_{\mbb C}^\us(G_{\mbb C},\chi_{_{\mbb C}})\geq 2$, and we 
denote $Y_{\chi,K}:=\bbV_K\invq_{\chi{\lower.2ex\hbox{$_{_K}$}}}G_K
=\Proj K[\bbV_K]^{G_K,\chi{\lower.2ex\hbox{$_{_K}$}}}$. 

Finally, we consider a character $\veps:\bG\rar\bG_m$ such that 
$\veps_{\mbb C}\in C(\chi_{\mbb C})$. Corollary \ref{cor:cone} 
implies that there is an open subset $U\subset\Spec\mbb Z$ such 
that $\bV^\sst(\bG,\chi)|_U\subset\bV^\sst(\bG,\veps)|_U$, and 
consequently we obtain a commutative diagram similar to \eqref{cd}, 
defined over $U$. 

\begin{proposition}{\label{thm:coh-line-bdls}}
There is a prime number $p_0\in\mbb N$ such that for any algebraically 
closed field $K$ of characteristic zero or ${\rm char\,}K>p_0$ holds: 

\nit{\rm (i)} $Y_{\chi,K}$ is arithmetically Cohen-Macaulay;
\smallskip

\nit{\rm (ii)}
$\begin{array}[t]{lll}
H^i\bigl(
Y_{\chi,K},
\psi_K^*\crl L_{Y_\veps,\veps^n{\kern-1.1ex\lower.25ex\hbox{$_{_K}$}}}
\bigr)=0
&
\forall\,n\in\mbb Z,
&
\forall\, i>0,\;i\neq\dim Y_{\veps, K};
\\[1ex] 
H^{\dim Y_{\veps,K}}\bigl(
Y_{\chi,K},
\psi_K^*\crl L_{Y_\veps,\veps^n{\kern-1.1ex\lower.25ex\hbox{$_{_K}$}}}
\bigr)=0
&
\forall\,n\geq 0.
&
\end{array}$
\end{proposition}

We recall from \cite[section 14]{hr} that a projective scheme is 
called arithmetically Cohen-Macaulay if it can be written as the 
`$\Proj$' of a graded, Cohen-Macaulay ring. 

\begin{proof}
For the beginning, we prove the result in characteristic zero, 
and we will conclude using a semi-continuity argument. 

So, let us fix an algebraically closed field $K$ of characteristic 
zero. Since in characteristic zero reductive groups are linearly 
reductive, the Hochster-Roberts theorem \cite[theorem 0.1]{ke1} 
implies that 
$$
R^{\bG_K}={(R^{\bG_K})}_0\oplus {(R^{\bG_K})}_+
:=K\oplus\bigoplus_{n>0}
K[\bbV_K]^{\bG_K}_{\chi^n{\kern-1.1ex\lower.25ex\hbox{$_{_K}$}}}
=K[\bbA^1_K\times\bbV_K]^{\bG_K}
$$
is Cohen-Macaulay, which proves the first statement. 

The proof of the second statement is done in two steps. We start 
proving it in the case when 
$\veps_{\mbb C}\in{\rm int.}\,C(\chi_{\mbb C})$, that is 
$C(\veps_{\mbb C})=C(\chi_{\mbb C})$. Corollary \ref{cor:cone} 
implies that the equality is valid still valid for $\veps_K$ and 
$\chi_K$, so that the quotients $Y_{\veps,K}$ and $Y_{\chi,K}$ 
coincide. Consequently we may assume that $\veps=\chi$. 
The Grothendieck-Serre correspondence \cite[theorem 20.4.4]{bs} 
gives the graded isomorphisms 
$$
\bigoplus_{n\in\mbb Z} 
H^i\bigl(Y_{\chi,K},\crl O_{Y_{\chi,K}}(n)\bigr)
\cong 
H^{\lower0.5ex\hbox{$^{i+1}$}}_{{(R^{\bG_K})}_+}(R^{\bG_K}),
\quad\forall\,i>0, 
$$
between the \v Cech and the local cohomology modules. On the 
other hand, since $\dim Y_{\chi,K}=\dim R^{\bG_K}-1$, we deduce 
from \cite[theorem 4.6]{ke1} that  
$$
H^{\lower0.5ex\hbox{$^{i}$}}_{{(R^{\bG_K})}_+}(R^{\bG_K})
=0,\text{ for }0<i\leq\dim Y_{\chi,K},
\quad\text{and}\quad
\biggl[
H^{\lower0.5ex\hbox
{$^{\dim Y_{\chi,K}+1}$}}_{{(R^{\bG_K})}_+}(R^{\bG_K})
\biggr]_{\geq 0}=0, 
$$
and our claim follows because 
$\crl O_{Y_{\chi,K}}(n)\cong
\crl L_{\chi^n{\kern-1.1ex\lower.25ex\hbox{$_{_K}$}}}$. 

Now we consider the general case, when 
$\veps_{\mbb C}\in C(\chi_{\mbb C})$; then, corollary \ref{cor:cone} 
implies that the same holds over $K$. We know from lemma \ref{lm:pre} 
that $\rR^i{\psi_K}_*\crl O_{Y_{\chi,K}}=0$, and we deduce  
$\rR^i{\psi_K}_*
\bigl(
\psi_K^*\crl L_{\veps^n{\kern-1.1ex\lower.25ex\hbox{$_{_K}$}}}
\bigr)
=
\crl L_{\veps^n{\kern-1.1ex\lower.25ex\hbox{$_{_K}$}}}
\otimes 
\rR^i{\psi_K}_*\crl O_{Y_{\chi,K}}
=0$. Therefore 
$$
H^i\bigl(Y_{\chi,K},
\psi_K^*\crl L_{\veps^n{\kern-1.1ex\lower.25ex\hbox{$_{_K}$}}}
\bigr)
\cong 
H^i\bigl(Y_{\veps,K},
\crl L_{\veps^n{\kern-1.1ex\lower.25ex\hbox{$_{_K}$}}}
\bigr),
$$
and our claim follows by applying the previous step. 

So far we have been concerned with fields of characteristic zero, 
and now we wish to extend the results in large characteristic. 
Since the morphism $\pr:\bY_\chi\!\rar\Spec\mbb Z$ is proper and 
flat, and $R^{\bG_{\mbb Q}}$ is Cohen-Macaulay,  
\cite[theorem 12.2.4 (i)]{ega4} implies that there is a non-empty 
open subset $U(\chi)\subset\Spec\mbb Z$ such that the homogeneous 
rings of the fibres over $U(\chi)$ are still Cohen-Macaulay; by 
flat base change \cite[lemma 2]{se}, $R^{\bG_K}$ is Cohen-Macaulay 
for all fields $K$ with $\lan{\rm char\,}K\ran\in U(\chi)$. This 
statement appears also in \cite[remark 1.3]{hr}. 

Before starting the proof of our second claim, let us remark that 
the difficulty relies in showing the {\em simultaneous} vanishing 
of the cohomology groups of the $\crl L_{\chi^n_K}$'s, for variable 
$n$. For fixed $n$, our claim is a direct application of the upper 
semi-continuity theorem \cite[theorem 12.8]{ha}. 

The proof of our second claim is divided again in two steps: if 
$\veps_{\mbb C}\in{\rm int.}\,C(\chi_{\mbb C})$, we may assume 
as before that $\veps=\chi$. Our previous discussion implies 
the vanishing 
$$
H^{\lower0.5ex\hbox{$^{i}$}}_{{(R^{\bG_K})}_+}(R^{\bG_K})
=0,\text{ for }0<i\leq\dim Y_{\chi,K},
$$
which, {\it via} the Grothendieck-Serre correspondence, 
translates into 
\begin{align}{\label{eqn:zero}}
H^i\bigl(Y_{\chi,K},
\crl L_{\chi^n{\kern-1.1ex\lower.25ex\hbox{$_{_K}$}}}
\bigr)
=0,\quad\forall\,0<i<\dim Y_{\chi,K},\;\forall\,n\in\mbb Z. 
\end{align}
It remains to prove the vanishing of the highest cohomology 
groups. Let us choose $c>0$ large enough such that 
$\mbb Z[\bV]^\bG_{\chi^{nc}}=
\bigl(\mbb Z[\bV]^\bG_{\chi^{c}}\bigr)^n$, and consider the 
corresponding projective embedding
$$
\xymatrix@C=0.3cm@R=0.5cm{
\bY_\chi
\ar[dr]_-{\pr}\ar@{^(->}[rr]^-{|\crl L_{\chi^c}|}
&& 
\mbb P^N_{\mbb Z}
\ar[dl]
\\ 
&\Spec\mbb Z&
}
$$
Since $\pr$ is proper and flat, the Hilbert polynomial of the 
ideal sheaves of the fibres $Y_{\chi,\mbb Z/\!\lan p\ran}=
\pr^{-1}\lan p\ran\hra\mbb P^N_{\mbb Z/\!\lan p\ran}$ is constant, 
so that there is a positive integer $m>0$ such that 
$\crl O_{Y_{\chi,\mbb Z/\!\lan p\ran}}$ is $m$-regular for all primes 
$p$ (see $e.g.$ \cite[page 101]{mu}). In particular 
\begin{align}{\label{eqn:max}}
H^{\dim Y_{\chi,\mbb Z/\!\lan p\ran}}
\bigl(
Y_{\chi,\mbb Z/\!\lan p\ran},
\crl L_{\chi^{nmc}_{\mbb Z/\!\lan p\ran}}
\bigr)=0,\quad\forall\,n>0. 
\end{align}
We consider now the diagonal representation 
$\bG\times_{\mbb Z}\bG\rar\bGl(M\oplus M)$, and consider the 
character $\chi\lower.35ex\hbox{$\times$}\chi^{mc}:
\bG\times_{\mbb Z}\bG\rar\bG_m$. The corresponding invariant 
ring is 
$$
R^{\bG\times\bG}
:=\mbb Z\oplus
\bigoplus_{n\geq 1}
\mbb Z[\bV]^\bG_{\chi^n}
\otimes_{\mbb Z}
\mbb Z[\bV]^\bG_{\chi^{nmc}}. 
$$
and we observe that $\Proj R^{\bG\times\bG}\cong 
\bY_\chi\times_{\mbb Z}\bY_\chi$, and that 
$\bigl(R^{\bG\times\bG}[n]\bigr)\!
\raise.1ex\hbox{$\widetilde{\phantom{xi}}$}\cong 
\crl L_{\chi^n}\boxtimes\crl L_{\chi^{nmc}}$. 
Applying the Hochster-Roberts theorem again, we deduce that 
$R^{\bG_{\mbb Q}\times\bG_{\mbb Q}}$ is Cohen-Macaulay, so 
that $R^{\bG_{K}\times\bG_{K}}$ is still Cohen-Macaulay for 
${\rm char\,}K$ large enough. Since $\dim Y_{\chi,K}< 
2\dim Y_{\chi,K}$ (when the quotient is a point there is nothing 
to prove), we deduce on the one hand that 
$$
H^{\dim Y_{\chi,K}}
\bigl(
Y_{\chi,K}\times Y_{\chi,K}, 
\crl L_{\chi^n}\boxtimes\crl L_{\chi^{nmc}}
\bigr)
=0,\quad\forall\,n>0
$$
On the other hand, equalities \eqref{eqn:zero} and \eqref{eqn:max} 
imply that 
$$
H^{\dim Y_{\chi,K}}\!
\bigl(
Y_{\chi,K}\times Y_{\chi,K}, 
\crl L_{\chi^n}\boxtimes\crl L_{\chi^{nmc}}\kern-.3ex
\bigr)
=
H^{\dim Y_{\chi,K}}\!
\bigl(Y_{\chi,K},\crl L_{\chi^n}\kern-.3ex\bigr)
\otimes_K
H^0
\bigl(Y_{\chi,K},\crl L_{\chi^{nmc}}\kern-.3ex\bigr).
$$ 
It follows that $H^{\dim Y_{\chi,K}}
\bigl(Y_{\chi,K},\crl L_{\chi^n}\bigr)=0$, for $n>0$. 
This finishes the proof of the second claim in the case 
$\veps_{\mbb C}\in{\rm int.\,}C(\chi_{\mbb C})$. 
 
For proving the general case, we are going to show that  
$\rR^i{\psi_K}_*\crl O_{Y_\chi,K}=0$, for all $i>0$, as soon as the 
characteristic of $K$ is large enough. The conclusion will follow 
then from the projection formula and the particular case treated 
above. 

Let us consider a finite set $\cal S\subset\bigoplus_{n\geq 1}
\bigl(\Sym^\bullet_{\mbb Z}M\bigr)^\bG_{\veps^n}$ of homogeneous 
elements (for this grading), which define an open, affine covering 
of $\bV^\sst(\bG,\chi)$. For an element $f\in\cal S$, we consider 
the diagram
$$
\xymatrix{
\psi_{\mbb C}^{-1}D(f_{\mbb C},Y_{\veps,\mbb C})
\ar[d]\ar[r]
&
\psi_{\mbb Q}^{-1}D(f_{\mbb Q},Y_{\veps,\mbb Q})
\ar[d]\ar[r]
&
\psi_U^{-1}D(f,\bY_{\veps})
\ar[d]
\\ 
\Spec\mbb C
\ar[r]
&
\Spec\mbb Q
\ar[r]
&
U
}
$$ 
We have proved in lemma \ref{lm:pre} that 
$H^i\bigl(\psi_{\mbb C}^{-1}D(f_{\mbb C},Y_{\veps,\mbb C}),
\crl O_{Y_\chi,\mbb C}\bigr)=0$, for all $i>0$. 
Applying \cite[proposition 9.3]{ha} to the faithfully flat 
base change $\Spec\mbb C\rar\Spec\mbb Q$, we deduce that 
$H^i\bigl(\psi_{\mbb Q}^{-1}D(f_{\mbb Q},Y_{\veps,\mbb Q}),
\crl O_{Y_\chi,\mbb Q}\bigr)=0$, for all $i>0$. Now we remark 
that $\crl O_{\bY_\chi}$ is flat over $\Spec\mbb Z$ (cf. remark 
\ref{rmk:hope}), so that the upper semi-continuity theorem 
\cite[theorem 12.8]{ha} applies: it implies that there 
is an open subset $U(f)\subset U\subset\Spec\mbb Z$ such that 
$H^i\bigl(\psi_{K}^{-1}D(f_{K},Y_{\veps,K}),
\crl O_{Y_\chi,K}\bigr)=0$, for all $i>0$, as soon as 
$\lan{\rm char}\,K\ran\in U(f)$. 
We let $U$ to be the finite intersection of the $U(f)$'s obtained 
this way. Applying \cite[proposition 9.3 and proposition 8.5]{ha}, 
we deduce the vanishing of 
$\rR^i{\psi_K}_*\crl O_{Y_\chi,K}|_{D(f_K,Y_{\veps,K})}$, 
for all $i>0$ and $f\in\cal S$, as soon as 
$\lan{\rm char}\,K\ran\in U$. 

Since $\bigl\{D(f_K,Y_{\veps,K})\bigr\}_{f\in\cal S}$ defines an 
open, affine covering of $Y_{\veps,K}$, we deduce that all the 
higher direct images $\rR^i{\psi_K}_*\crl O_{Y_\chi,K}$, $i>0$, 
vanish.  
\end{proof}

\begin{remark}{\label{rmk:fail}}
This proposition gives very satisfactory results about the vanishing 
of the higher cohomology groups of ample line bundles on our geometric 
quotients, but gives only partial answer in the case of nef line 
bundles. More precisely, as an immediate consequence, we deduce 
that for $\veps\in\partial C(\chi)\cap\cal X^*(G)$ there is a 
constant $c>0$ (depending on $\veps$) such that 
$$
H^i(Y_\chi,\crl L_{Y_\chi,\veps^{c}}^{\otimes n})=0,
\; 
\forall n\geq 0,
\quad\text{and}\quad
H^i(Y_\chi,\crl L_{Y_\chi,\veps^{-c}}^{\otimes n})=0,
\; 
\forall n\geq 0,\;i\neq\dim Y_\chi. 
$$
Indeed, if we choose $c>0$ is such that 
$\mbb Q[\bbV_{\mbb Q}]^{G_{\mbb Q},\veps_{\mbb Q}^c}$ is generated 
in degree one, then $\crl L_{Y_\veps,\veps_{\mbb Q}^c}$ is locally 
free, so that $\psi^*_{\mbb Q}\crl L_{Y_\veps,\veps_{\mbb Q}^c}=
\crl L_{Y_\chi,\veps_{\mbb Q}^c}$ by lemma \ref{glob-gen}. 

The trouble is that this statement is unnatural, since we would 
like to have the cohomology vanishing for all powers, not only 
for the multiples of some number. This is what we are going to 
prove next. 
\end{remark}

So, let us assume that the character $\veps$ has the property that 
$\crl L_{Y_\chi,\veps_{\mbb C}}\rar Y_{\chi,\mbb C}$ is invertible, and 
$\veps_{\mbb C}\in\partial C(\chi_{\mbb C})\cap\cal X^*(G_{\mbb C})$. 
Then the same holds (by faithfully flat base change) for 
$\crl L_{Y_\chi,\veps_{\mbb Q}}\rar Y_{\chi,\mbb Q}$, which implies 
that $\crl L_{\bY_\chi,\veps}|_{U'(\veps)}\rar\bY_\chi|_{U'(\veps)}$ 
is still invertible over some non-empty, open subset 
$U'(\veps)\subset\Spec\mbb Z$. 

We define the representation 
$$
\bar\rho:\bar\bG:=\bG\times_{\Spec\mbb Z}\bG_m
\lar
\bGl\bigl(M\oplus\mbb Z^{\oplus 2}\bigr),
\quad 
\bar\rho(g,t):={\rm diag}\bigl(\rho(g),t\veps^{-1}(g),t\bigr),
$$
and notice that since $\bigl(\Sym^\bullet_{\mbb Z}M\bigr)^\bG=\mbb Z$, 
$\bigl(\Sym^\bullet_{\mbb Z}(M\oplus\mbb Z^2)\bigr)^{\bar\bG}
=\mbb Z$ too. 

After replacing $\chi$ with a suitably large positive multiple, 
we may and we assume that 
$\bigl(\Sym^\bullet_{\mbb Z}M\bigr)^{\bG,\chi}$ is generated in 
degree one. Let us define now the character 
$$
\bar\chi:\bar\bG\lar\bG_m,\quad\bar\chi(g,t):=\chi(g)t, 
$$
and observe that 
\begin{align}{\label{eqn:2}}
\begin{array}{ll}
\bigl(\Sym^\bullet_{\mbb Z}(M\oplus\mbb Z^2)\bigr)^{\bar\bG,\bar\chi}
&
=
\bigl((\Sym^\bullet_{\mbb Z}M)[w_1,w_2]\bigr)^{\bar\bG,\bar\chi} 
= 
\underset{n\geq 0}{\bigoplus}
\bigl((\Sym^\bullet_{\mbb Z}M)[w_1,w_2]\bigr)^{\bar\bG}_{\bar\chi^n} 
\\[2ex] 
&
=
\underset{n\geq 0}{\bigoplus}\;
\underset{a+b=n}{\bigoplus}
(\Sym^\bullet_{\mbb Z}M)^\bG_{\chi^n\veps^a}w_1^aw_2^b. 
\end{array}
\end{align}
We have proved in corollary \ref{cor:cone} that there is an open 
subset $U''(\veps)\subset\Spec\mbb Z$ such that $\veps_K\in C(\chi_K)$ 
for any algebraically closed field $K$ having 
$\lan{\rm char\,}K\ran\in U''(\veps)$. We define 
$U(\veps):=U'(\veps)\cap U''(\veps)$; it is a non-empty, open subset 
of $\Spec\mbb Z$. More precisely, $U(\veps)=\Spec B$, where $B$ is 
obtained out of $\mbb Z$ by inverting finitely many primes.
 
\begin{lemma}{\label{lm:p1}}
Let $\veps:\bG\rar\bG_m$ be a character such that 
$\veps_{\mbb C}\in C(\chi_{\mbb C})\cap\cal X^*(G_{\mbb C})$, 
and moreover $\crl L_{Y_\chi,\veps_{\mbb C}}\rar Y_{\chi,\mbb C}$ 
is invertible. Then there is a non-empty, open subset 
$U(\veps)\subset\Spec\mbb Z$ over which the following isomorphism 
holds 
$$
\mbb P\bigl(
\crl O_{\bY_\chi}\oplus\crl L_{\bY_\chi,\veps}
\bigr)|_{U(\veps)}
:=
\uProj\bigl(
\Sym^\bullet_{\crl O_{\bY_\chi}}
(\crl O_{\bY_\chi}\oplus\crl L_{\bY_\chi,\veps})
\bigr)
\cong 
\Proj\bigl(
\mbb Z[\bV\times_{\mbb Z}\mbb A^2_{\mbb Z}]^{\bar\bG,\bar\chi}
\bigr)|_{U(\veps)}.
$$
\end{lemma}

\begin{proof} 
We prove the statement for $U(\veps)=\Spec B$ defined above. For shorthand, 
we write $\crl O:=\crl O_{\bY_\chi}$, $\crl L:=\crl L_{\bY_\chi,\veps}$, 
and $B[\bV]:=\Sym^\bullet_{B}(M\otimes_{\mbb Z}B)$. 
Since we have assumed that $\mbb Z[\bV]^{\bG,\chi}$ is generated in degree 
one, there is a finite subset $\{f_j\}_{j\in J}\subset\mbb Z[\bV]^\bG_\chi$ 
such that the non-vanishing loci $\{D(f_j,\bY_\chi)\}_{j\in J}$ 
form an open, affine covering of $\bY_\chi$. We notice that 
$D(f_j,\bY_\chi|_{U(\veps)})=\Spec\bigl((f_j^{-1}B[\bV])^\bG\bigr)$. By 
our previous discussion, $\crl L|_{U(\veps)}\rar\bY_\chi|_{U(\veps)}$ is 
invertible, so that we find trivializing sections 
$s_j\in\bigl(f_j^{-1}B[\bV]\bigr)^\bG_\veps$, $j\in J$, such that 
$$
\crl L|_{D(f_j,\bY_\chi|_{U(\veps)})}=
\bigl(s_j \cdot( f_j^{-1}B[\bV])^\bG\bigr)\!
\raise.1ex\hbox{$\widetilde{\phantom{iii}}$}.
$$ 
Choosing $f$ to be one of the $f_j$'s and $s$ the corresponding $s_j$, 
we find that 
\begin{align}{\label{eqn:1}}
\begin{array}{ll}
\mbb P\bigl(
\crl O_{\bY_\chi}\oplus\crl L_{\bY_\chi,\veps}
\bigr)|_{D(f_j,\bY_\chi|_{U(\veps)})}
&
\\[1ex]
&
\kern-12em
=\Proj
\biggl(
\Sym^\bullet_{(f^{-1}B[\bV])^\bG}
\bigl(
\bone\cdot(f^{-1}B[\bV])^\bG\oplus s_j \cdot( f_j^{-1}B[\bV])^\bG
\bigr)
\biggr)
\\[1.5ex]
&
\kern-12em
=\Proj
\biggl(
\underset{n\geq 0}{\bigoplus}\;
\underset{a+b=n}{\bigoplus}
(f^{-1}B[\bV])^\bG\cdot s^a\bone^b
\biggr)
=\Proj
\biggl(
(f^{-1}B[\bV])^\bG\bigl[s,\bone\bigr]
\biggr).
\end{array}
\end{align}
On the other hand, we know that $\veps_K\in C(\chi_K)$ for all 
algebraically closed field with $\lan{\rm char\,}K\ran\in U(\veps)$. 
It follows that for all $n\geq 1$ and $0\leq a\leq n$, 
$\chi^n_K\veps^a_K\in{\rm int.\,}C(\chi_K)$. We deduce now from 
\eqref{eqn:2} that the $\bar\chi$-semi-stable locus of 
$(\bV\times_{\mbb Z}\mbb A^2_{\mbb Z})|_{U(\veps)}$ is 
$$
{\bigl(
(\bV\times_{\mbb Z}\mbb A^2_{\mbb Z})|_{U(\veps)}
\bigr)}^\sst(\bar\bG,\bar\chi)
=
{\bigl(
\bV|_{U(\veps)}
\bigr)}^\sst(\bG,\chi)
\times_{\mbb Z}
(\mbb A^2_{\mbb Z}\sm\{0\}). 
$$ 
In particular, the $\bar\chi$-unstable locus has codimension two; 
since both $\crl L_{\bY_\chi,\chi}|_{U(\veps)}$ and 
$\crl L_{\bY_\chi,\veps}|_{U(\veps)}$ are invertible, remark 
\ref{rmk:inj} implies that 
\begin{align}{\label{eqn:3}}
\crl L_{\bY_\chi,\chi}^{\otimes n}|_{U(\veps)}
\otimes 
\crl L_{\bY_\chi,\veps}^{\otimes a}|_{U(\veps)}
\cong 
\crl L_{\bY_\chi,\chi^n\veps^a}|_{U(\veps)},\quad 
\forall\,n,a\in\mbb Z. 
\end{align}
The categorical quotient of this scheme for the $\bar\bG$-action is 
obtained by glueing the categorical quotients of the $\bar\bG$-invariant 
open subschemes 
$D(f_j)\times_{\mbb Z}(\mbb A^2_{\mbb Z}\sm\{0\})|_{U(\veps)}$. 
Again, taking $f$ to be one of the $f_j$'s, we observe that 
\begin{align}{\label{eqn:4}}
\begin{array}{r}
\kern-.3ex\Proj\biggl(
\underset{n\geq 0}{\bigoplus}\;
\underset{a+b=n}{\bigoplus}\kern-.9ex
(f^{-1}B[\bV])^\bG_{\chi^n\veps^a}w_1^aw_2^b
\biggr)
\kern-.9ex\srel{\eqref{eqn:3}}{=}\kern-.5ex
\Proj\biggl(
\underset{n\geq 0}{\bigoplus}\;
\underset{a+b=n}{\bigoplus}\kern-.85ex
(f^{-1}B[\bV])^\bG f^n s^a w_1^a w_2^b
\biggr)
\\[1.5ex]
=
\Proj\biggl(
\underset{n\geq 0}{\bigoplus}\;
\underset{a+b=n}{\bigoplus}\kern-.9ex
(f^{-1}B[\bV])^\bG \bigl[fsw_1,fw_2\bigr]
\biggr). 
\end{array}
\end{align}
The equalities \eqref{eqn:1} and \eqref{eqn:4} show that the varieties 
we are considering are locally isomorphic, and one can also check that 
these local isomorphisms are compatible with coordinate changes. 
\end{proof}

Now we are in position to prove our main cohomology vanishing result. 

\begin{theorem}{\label{thm:zero}}
Assume that our previous notations are in force. Let $\bG\rar\bGl(M)$ 
be a representation such that $(\Sym^\bullet_{\mbb Z}M)^\bG=\mbb Z$, 
and consider two characters $\chi,\veps$ of $\bG$ with the properties: 

-- $\codim_{\bbV_{\mbb C}}\bbV^\us(G_{\mbb C},\chi_{\mbb C})\geq 2$;

-- $\crl L_{Y_\chi,\veps_\mbb C}\lar Y_{\chi,\mbb C}$ is invertible 
and nef. 

\nit Then there is a prime $p_0$ such that for all algebraically closed 
fields of characteristic ${\rm char\,}K>p_0$ holds: 

$\begin{array}{lll}  
\hbox{-- } H^i\bigl(
Y_{\chi,K},{\crl L}^{\otimes n}_{Y_\chi,\veps_K}
\bigr)=0
&
\forall\,n\geq 0,
&
\forall\, i>0;
\\[1ex] 
\hbox{-- } H^i\bigl(
Y_{\chi,K},{\crl L}^{\otimes n}_{Y_\chi,\veps_K^{-1}}
\bigr)=0
&
\forall\,n>0,
&
\forall\, i\neq\dim Y_{\veps, K}
\end{array}
$
\end{theorem}

\begin{proof}
Let $U(\veps)\subset\Spec\mbb Z$ be as in the previous lemma, and 
write $\crl L:=\crl L_{Y_\chi,\veps}$. Then 
$\crl L|_{U(\veps)}\rar\bY_\chi|_{U(\veps)}$ is relatively nef, in 
the sense that is nef on the geometric fibres of 
$\bY_\chi|_{U(\veps)}\rar U(\veps)$. We claim that 
$$
\crl O_{\mbb P(\crl O\oplus\crl L)}(1)|_{U(\veps)}
\lar 
\mbb P\bigl(\crl O\oplus\crl L\bigr)|_{U(\veps)}
$$
is still relatively nef. Indeed, let us consider a relatively ample, 
invertible sheaf $\crl A\rar\bY_\chi$. Then for any integer $k>0$ 
holds 
$$
\xymatrix@C=0.1cm@R=0.5cm{
\crl A\otimes\crl O_{\mbb P(\crl O\oplus\crl L)}(k)
=\jmath_k^*\kern1.5ex\tld{\kern-1.5ex\crl A}_k
\ar[d]
&&
\tld{\kern-1.5ex\crl A}_k:=\crl O_{\mbb P\left(\crl A\otimes
(\crl O\oplus\crl L\oplus{\dots}\oplus\crl L^k)\right)}(1)
\ar[d]
\\
\mbb P(\crl O\oplus\crl L)\;
\ar@{^(->}[rr]^-{\jmath_k}\ar[dr]
&&
\mbb P\bigl(\crl A\otimes
(\crl O\oplus\crl L\oplus{\dots}\oplus\crl L^k)\bigr)
\ar[dl]
\\
&\bY_\chi|_{U(\veps)}&
}
$$
According to \cite[proposition 2.2 and 3.2]{ha1}, 
$\kern1.5ex\tld{\kern-1.5ex\crl A}_k$ is relatively ample. As $k>0$ 
is arbitrary, we deduce that $\crl O_{\mbb P(\crl O\oplus\crl L)}(1)$ 
is in the closure of the relatively ample cone, and therefore is 
relatively nef. 

Proposition \ref{thm:coh-line-bdls}, together 
with lemma \ref{lm:p1}, imply that there is a constant $c>0$ having 
the property that, after possibly shrinking $U(\veps)$ further, the 
following vanishing holds 
$$
\begin{array}{lll}
H^i\bigl(
\mbb P\bigl(\crl O\oplus\crl L\bigr)_K,{\crl O(1)}^{nc}_{_K}
\bigr)=0
&
\forall\,n\geq 0,
\\[1,5ex]
H^i\bigl(
\mbb P\bigl(\crl O\oplus\crl L\bigr)_K,\crl O(-1)^{nc}_{_K}
\bigr)=0
&
\forall\,n\geq 0\text{ and }i\neq\dim Y_{\veps,K}+1, 
\end{array}
$$
on all geometric fibres over $U(\veps)$. Denoting  
$\pi:\mbb P\bigl(\crl O\oplus\crl L\bigr)\rar \bY_\chi$ the 
natural projection, we have  
$$
\rR^i{(\pi_K)}_*{\crl O(1)}^{nc}_{_K}=0,
\quad
\forall\,i>0\;\forall\,n\geq 0,
$$
and therefore 
$$
0=H^i\bigl(
\mbb P\bigl(\crl O\oplus\crl L\bigr)_K,{\crl O(1)}^{nc}_{_K}
\bigr)
=
H^i\bigl(
Y_{\chi,K},{(\pi_K)}_*{\crl O(1)}^{nc}_{_K}
\bigr)
=
\bigoplus_{a=0}^{nc}
H^i\bigl(
Y_{\chi,K},{\crl L}^{a}_{_K}
\bigr). 
$$
Taking $n>0$ arbitrarily large, we deduce the vanishing of the 
positive powers of nef bundles. 

For the negative powers we proceed as follows: the relative 
canonical sheaf of $\pi$ is 
$\omega_{\rm rel}\cong \pi^*\crl L\otimes \crl O(-2)$, so that 
the relative duality implies 
$$
\rR^1\pi_*{\crl O(-1)}^{nc}
\cong 
\pi_*\bigl(\pi^*\crl L\otimes\crl O(nc-2)\bigr)^\vee
=\bigoplus_{a=1}^{nc-1}{\crl L}^{-a}. 
$$
Since ${(\pi_K)}_*{\crl O(-1)}^{nc}=0$ for $n>0$, the Leray 
spectral sequence implies that 
$$ 
\begin{array}{rl}
H^{i+1}\bigl(
\mbb P\bigl(
\crl O\oplus\crl L\bigr)_K,{\crl O(-1)}^{nc}_{_K}
\bigr)
&
=
H^i\bigl(
Y_{\chi,K},\rR^1{(\pi_K)}_*{\crl O(-1)}^{nc}_{_K}
\bigr)
\\[1ex]
&\disp
=
\bigoplus_{a=1}^{nc-1}
H^i\bigl(
Y_{\chi,K},{\crl L}^{-a}_{_K}
\bigr).
\end{array}
$$
The left-hand-side vanishes for $i\neq\dim Y_{\veps,K}$. Taking 
again $n>0$ arbitrarily large, we deduce the cohomology vanishing 
for the negative powers. 
\end{proof}

\begin{remark}
It is actually possible `to squeeze out' some more information, namely 
that the prime $p_0$ appearing in the previous theorem can be chosen 
independent of $\veps$. In other words, one is able to find a prime 
$p(\chi)$, depending only on the quotient $\bY:=\bY_{\chi}$, for which 
one has vanishing of the higher cohomology groups, for all invertible, 
nef sheaves, in characteristic larger than $p(\chi)$. 

This can be seen as follows: consider a {\em finite} set of characters 
$\{\veps_j\}_{j\in J}$ such that $\{\crl L_{\veps_j}\}_{j\in J}$ 
generate the nef cone $C(\chi_{\mbb Q})$ of $\bY_{\chi,\mbb Q}$ over 
$\mbb Z_{\geq 0}$; then the same holds over a non-empty, open subset 
of $\Spec\mbb Z$. Similarly as in lemma \ref{lm:p1}, one proves that 
$$
\hbox{
$\mbb P\bigl(\crl O\oplus\underset{j\in J}{\bigoplus}\crl L_j\bigr)$
}
\cong 
\bigl[\bV\times_{\mbb Z}\underbrace{\mbb A_{\mbb Z}^2\times_{\mbb Z}
{\dots}
\times_{\mbb Z}\mbb A_{\mbb Z}^2}_{|J|\text{ times}}\,\bigr]
\bigr/\kern-.8ex\bigr/
\bG\times_{\mbb Z}\bG_m^{|J|}
$$ 
for a suitable action and linearization. Applying now proposition 
\ref{thm:coh-line-bdls} to large multiples of 
$\crl O_{\mbb P(\crl O\oplus\underset{j\in J}{\bigoplus}\crl L_j)}(1)$, 
just as in \ref{thm:zero}, one obtains the desired result.  
\end{remark}

We have focused on invariant quotients of affine spaces for 
actions of reductive groups. However, we observe that the 
Grothendieck-Serre correspondence, together with 
\cite[theorem 0.1 and theorem 3.4]{ke1} immediately imply the 

\begin{theorem}{\label{thm:X}}
Let $X$ be a smooth, affine variety, defined over an algebraically 
closed field $K$ of characteristic zero, which is acted on by a 
reductive group $G$. We assume that the $K[X]^G=K$, and consider 
a character $\chi\in\cal X^*(G)$ such that 
$X^\sst(G,\chi)=X^\s_{(0)}(G,\chi)$ and 
${\rm codim\,}_X X^\us(G,\chi)\geq 2$. 
Then: 

\nit{\rm (i)} $X\invq_\chi G$ is arithmetically Cohen-Macaulay; 

\nit{\rm (ii)} For an invertible and nef sheaf $\crl L\rar X$, 

$\begin{array}[t]{lll}
H^i\bigl(X\invq_\chi G,\crl L\bigr)=0
&
\forall\, i>0\,;
\\[1ex] 
H^i\bigl(X\invq_\chi G,\crl L^{-1}\bigr)=0
&
\forall\,i\neq\kappa(\crl L),
\end{array}$

\nit where $\kappa(\crl L)$ denotes the Kodaira-Iitaka dimension of 
$\crl L$. 
\end{theorem}

\begin{proof}
The first statement is just the usual Hochster-Roberts theorem.  
For our second claim, we start by noticing that the very same argument 
as in remark \ref{rmk:inj} shows that the pull-back by the quotient 
map induces the monomorphism
$$
0\lar\Pic(X\invq_\chi G)\lar\cal X^*(G). 
$$
Therefore $\crl L=\crl L_\veps$ for some $\veps\in\cal X^*(G)$ (the 
latter is defined similarly as in \ref{line-bdls}). The GIT-cone 
theorem \ref{prop:fan} implies that $\veps\in C(\chi)$, so that a 
positive multiple of $\crl L$ is globally generated. 

If $\veps\in{\rm int.\,}C(\chi)$, we may assume $\veps=\chi$, and 
the vanishings follow directly from the fact that the invariant 
ring $K[X]^{G,\chi}$ is Cohen-Macaulay. Otherwise we consider the 
quotient map $\psi:X\invq_\chi G\rar X\invq_\veps G$; there is $c>0$ 
and an ample line bundle $\crl L'\rar X\invq_\veps G$ such that 
$\psi^*\crl L'\cong\crl L^c$. Using an appropriate version of lemma 
\ref{lm:pre} and the previous step, we deduce that the cohomology 
vanishings hold for arbitrary multiples (positive and negative) 
of $\crl L^c$. 

In order to conclude, we use the same trick as before: namely, 
$\mbb P\bigl(\crl O\oplus\crl L\bigr)$ is the quotient of 
$X\times\mbb A^2_K$ for a suitable $G\times G_m$-action, and we 
may apply the previous step to the nef line bundle 
$\crl O_{\mbb P(\crl O\oplus\crl L)}(1)$. 
\end{proof}


\end{document}